\documentclass[12pt,leqno,draft]{article}
\usepackage{amsfonts}
\usepackage{amsmath, dsfont, amsthm, amsfonts, amssymb, color}
\usepackage{mathrsfs}
\usepackage{color}
\usepackage{stmaryrd}
\newtheorem{thm}{Theorem}[section]

\newtheorem{lem}[thm]{Lemma}

\newtheorem{exa}[thm]{Example}
\newtheorem{rem}[thm]{Remark}
\theoremstyle{definition}

\newcommand{\scr}[1]{\mathscr #1}
\definecolor{wco}{rgb}{0.5,0.2,0.3}

\numberwithin{equation}{section} \theoremstyle{remark}

\newcommand{\ua}{\uparrow}

\title{{\bf %Zvonkin's transform on $\R^d\times \scr P_2$ and Applications for
McKean-Vlasov SDEs with Bounded Measurable Interaction}\footnote{Supported in
 part by  National Key R\&D Program of China (2022YFA1006000), NNSFC (12271398, 12531007).} }
\author{
{\bf   Xing Huang  }\\
\footnotesize{ Center for Applied Mathematics, Tianjin
University, Tianjin 300072, China}\\
\footnotesize{  xinghuang@tju.edu.cn}}
\begin{document}
\allowdisplaybreaks
\def\R{\mathbb R}  \def\ff{\frac} \def\ss{\sqrt} \def\B{\mathbf
B} \def\W{\mathbb W}
\def\N{\mathbb N} \def\kk{\kappa} \def\m{{\bf m}}
\def\ee{\varepsilon}\def\ddd{D^*}
\def\dd{\delta} \def\DD{\Delta} \def\vv{\varepsilon} \def\rr{\rho}
\def\<{\langle} \def\>{\rangle} \def\GG{\Gamma} \def\gg{\gamma}
  \def\nn{\nabla} \def\pp{\partial} \def\E{\mathbb E}
\def\d{\text{\rm{d}}} \def\bb{\beta} \def\aa{\alpha} \def\D{\scr D}
  \def\si{\sigma} \def\ess{\text{\rm{ess}}}
\def\beg{\begin} \def\beq{\begin{equation}}  \def\F{\scr F}
\def\Ric{\text{\rm{Ric}}} \def\Hess{\text{\rm{Hess}}}
\def\e{\text{\rm{e}}} \def\ua{\underline a} \def\OO{\Omega}  \def\oo{\omega}
 \def\tt{\tilde} \def\Ric{\text{\rm{Ric}}}
\def\cut{\text{\rm{cut}}} \def\P{\mathbb P} \def\ifn{I_n(f^{\bigotimes n})}
\def\C{\scr C}      \def\aaa{\mathbf{r}}     \def\r{r}
\def\gap{\text{\rm{gap}}} \def\prr{\pi_{{\bf m},\varrho}}  \def\r{\mathbf r}
\def\Z{\mathbb Z} \def\vrr{\varrho} \def\ll{\lambda}
\def\L{\scr L}\def\Tt{\tt} \def\TT{\tt}\def\II{\mathbb I}
\def\i{{\rm in}}\def\Sect{{\rm Sect}}  \def\H{\mathbb H}
\def\M{\scr M}\def\Q{\mathbb Q} \def\texto{\text{o}} \def\LL{\Lambda}
\def\Rank{{\rm Rank}} \def\B{\scr B} \def\i{{\rm i}} \def\HR{\hat{\R}^d}
\def\to{\rightarrow}\def\l{\ell}\def\iint{\int}
\def\EE{\scr E}\def\Cut{{\rm Cut}}
\def\A{\scr A} \def\Lip{{\rm Lip}}
\def\BB{\scr B}\def\Ent{{\rm Ent}}\def\L{\scr L}
\def\R{\mathbb R}  \def\ff{\frac} \def\ss{\sqrt} \def\B{\mathbf
B}
\def\N{\mathbb N} \def\kk{\kappa} \def\m{{\bf m}}
\def\dd{\delta} \def\DD{\Delta} \def\vv{\varepsilon} \def\rr{\rho}
\def\<{\langle} \def\>{\rangle} \def\GG{\Gamma} \def\gg{\gamma}
  \def\nn{\nabla} \def\pp{\partial} \def\E{\mathbb E}
\def\d{\text{\rm{d}}} \def\bb{\beta} \def\aa{\alpha} \def\D{\scr D}
  \def\si{\sigma} \def\ess{\text{\rm{ess}}}
\def\beg{\begin} \def\beq{\begin{equation}}  \def\F{\scr F}
\def\Ric{\text{\rm{Ric}}} \def\Hess{\text{\rm{Hess}}}
\def\e{\text{\rm{e}}} \def\ua{\underline a} \def\OO{\Omega}  \def\oo{\omega}
 \def\tt{\tilde} \def\Ric{\text{\rm{Ric}}}
\def\cut{\text{\rm{cut}}} \def\P{\mathbb P} \def\ifn{I_n(f^{\bigotimes n})}
\def\C{\scr C}      \def\aaa{\mathbf{r}}     \def\r{r}
\def\gap{\text{\rm{gap}}} \def\prr{\pi_{{\bf m},\varrho}}  \def\r{\mathbf r}
\def\Z{\mathbb Z} \def\vrr{\varrho} \def\ll{\lambda}
\def\L{\scr L}\def\Tt{\tt} \def\TT{\tt}\def\II{\mathbb I}
\def\i{{\rm in}}\def\Sect{{\rm Sect}}  \def\H{\mathbb H}
\def\M{\scr M}\def\Q{\mathbb Q} \def\texto{\text{o}} \def\LL{\Lambda}
\def\Rank{{\rm Rank}} \def\B{\scr B} \def\i{{\rm i}} \def\HR{\hat{\R}^d}
\def\to{\rightarrow}\def\l{\ell}\def\BB{\mathbb B}
\def\8{\infty}\def\I{1}\def\U{\scr U} \def\n{{\mathbf n}}
\maketitle

\begin{abstract} In this paper, McKean-Vlasov SDEs with bounded measurable interaction is investigated. The regularity estimate
$$\|P_t^\ast\gamma^1-P_t^\ast\gamma^2\|_{var}\leq ct^{-\frac{1}{2}}\W_{1}(\gamma^1,\gamma^2),\ \ t\in(0,T]$$
for the nonlinear semigroup $P_t^\ast$ associated to McKean-Vlasov SDEs is derived. Two cases are considered respectively.
The first case concentrates on the model where the interaction in the drift is merely assumed to be bounded measurable while the distribution dependent diffusion term is allowed to be Lipschitz continuous under $L^\eta$($\eta\in(0,1)$)-Wasserstein distance in the measure variable. In the second case, the diffusion is distribution free and the drift contain two parts: a bounded measurable interaction term plus a partially dissipative term. As an application of the regularity estimate, the exponential ergodicity in $\W_1$ is obtained in the second case.
\end{abstract} \noindent
 \vskip 2cm

\section{Introduction}
\subsection{Backgrounds and motivation}

Let $\scr P$ be the set of all probability measures on $\R^d$ equipped with the weak topology, and $W_t$ be an $m$-dimensional Brownian motion on a complete filtration probability space $(\OO,\{\F_t\}_{t\ge 0},\F,\P)$. Consider
the following McKean-Vlasov SDE arising from the pioneering work in \cite{McKean} on $\R^d$:
\beq\label{E0} \d X_t= b_t(X_t, \L_{X_t})\d t+  \si_t(X_t,\L_{X_t})\d W_t,\ \ t\geq 0,\end{equation}
where   $\L_{X_t}$ is the distribution of $X_t$, and
$$b: [0,\infty)\times\R^d\times\scr P\to\R^d,\ \ \si: [0,\infty)\times\R^d\times \scr P\to \R^d\otimes\R^m$$
are measurable. It is well known that the law of \eqref{E0} solves a nonlinear Fokker-Planck equation. Due to its wide applications, \eqref{E0} has been intensively investigated, see for instance \cite{BBP, CR,CF,FZ,GL, HWJMAA, L, MV,Ren, FYW1} and references therein for the well-posedess under various assumptions, \cite{B,CF,CM,HWJMAA,HW22b,Ren,S,FYW3} for the regularity including the derivative formula(estimate) and Wang's Harnack inequality and \cite{JW2016, L, SY, SZ} for the propagation of chaos and limit theorem.

We should remark that most of the aforementioned results concentrate on the case that the diffusion coefficient is distribution free. When the diffusion coefficient is distribution dependent, the diffusion coefficient changes once initial law differs,  so that the Girsanov transform is unavailable and the well-posedness as well as the regularity estimate becomes more difficult.

To overcome the difficulty caused by the distribution dependent diffusion, \cite{CR} adopted the method of parametrix expansion. \cite{CF} extends the result in \cite{CR} by introducing the linear functional derivative with respect to the measure argument. Also by the method of parametrix expansion, in \cite{HWJMAA}, the author and his collaborators proved the well-posedness of \eqref{E0} under the assumption that the diffusion is Lipschitz continuous under $L^k(k\geq 1)$-Wasserstein distance and the drift is Lipschitz continuous under the weighted variation distance plus $L^k(k\geq 1)$-Wasserstein distance. However, the method of parametrix expansion seems a bit complicated since a series is involved in.

Recently, in \cite[Theorem 1.3(1)]{HRWJDE}, the author and his collaborators established the estimate of probability distances for two different diffusion processes and apply it to investigate the well-posedness of \eqref{E0}, where the drift and the diffusion terms are Lipschitz continuous in the measure variable under the probability distance $\W_\psi$ defined by
$$\W_\psi(\mu,\nu):=\sup_{f:\ \ \sup_{x\ne y} \ff{|f(x)-f(y)|}{\psi(|x-y|)}\le 1} \left|\int_{\R^d}f\d \mu-\int_{\R^d}f\d \nu\right|$$
for some increasing and concave $\psi$ with $\lim_{t\to0}\psi(t)=0$.
Moreover, in \cite[Theorem 1.3(2)]{HRWJDE}, the authors gained the regularity estimate
$$\W_\psi(P_t^*\gg,P_t^*\tt\gg)\le \ff{c\psi(t^{\ff 1 2})}{\ss t} \W_k(\gg,\tt\gg),\ \ t\in(0,T],$$
where $P_t^*\gg$ is the distribution of the solution to \eqref{E0} with initial law $\gamma$ and $\psi$ satisfies the Dini condition
\beq\label{ESC} \int_0^1   \ff{ \psi(s)} s\,   \d s<\infty. \end{equation}
However, \cite[Theorem 1.3]{HRWJDE} excludes the critical case $\psi=2$, in which $\W_\psi$ coincides with the total variation distance $\|\cdot\|_{var}$.

In this paper, we aim to replace
$$|b_t(x,\gamma)-b_t(x,\tilde{\gamma})|\leq K\W_{\psi}(\gamma,\tilde{\gamma})$$
%for some increasing and concave function $\psi$ with $\lim_{t\to 0}\psi(t)=0$ in \cite[Theorem 1.3(2)]{HRWJDE}
with
$$|b_t(x,\gamma)-b_t(x,\tilde{\gamma})|\leq K\|\gamma-\tilde{\gamma}\|_{var}.$$
One may also expect to do this replacement for $\sigma$. However,
the following counterexample indicates that the uniqueness of \eqref{E0} may be destroyed if $\sigma$ is merely Lipschitz continuous in measure variable under total variation distance. Consider
\begin{align}\label{tyu}
\d X_t=\P(X_t\in\{0\}^{c})\d W_t.
\end{align}
In this case, $\sigma(\mu)=\mu(\{0\}^{c})$ and it is Lipschitz continuous in measure variable under total variation distance. Note that \eqref{tyu} has two solutions $X_t^1=0, X_t^2=W_t$ if $X_0=0$. Inspired by this counterexample, we will assume that $\sigma$ is Lipschitz continuous under $\W_\eta$ for some $\eta\in(0,1)$, see Assumption (A3) below for details.

Besides the well-posedness, we also focus more on the regularity estimate
\begin{align}\label{koz}\|P_t^\ast\gamma^1-P_t^\ast\gamma^2\|_{var}\leq ct^{-\frac{1}{2}}\W_{1}(\gamma^1,\gamma^2),\ \ t\in(0,T].
\end{align}
It plays an important role in the studying of exponential ergodicity for McKean-Vlasov SDEs with bounded measurable interaction.

To study the ergodicity, we consider the following time-homogeneous McKean-Vlasov SDEs on $\R^d$:
\beq\label{E0i} \d X_t= b (X_t, \L_{X_t})\d t+  \si(X_t)\d W_t,\ \ t\geq 0,\end{equation}
where
$$ \si: \R^d\to \R^d\otimes\R^m,\ \ \ b:  \R^d\times \scr P \to \R^d $$
are measurable. Assume $b(x,\mu)=b^{(1)} (x)+b^{(0)} (x, \mu)$.
When $b^{(1)} (x)$ is partially dissipative and $b^{(0)}(x,\mu)$ is Lipschitz continuous in $\mu$ under $L^1$-Wasserstein distance(for instance, the interaction kernel is Lipschitz continuous), the exponential  ergodicity of \eqref{E0i} in $L^1$-Wasserstein distance has been intensively investigated in \cite{EGZ,Schuh 2024, WFY 2023,ZSQ} and the references therein. The main tool is (asymptotic) reflecting coupling. Quite recently, in \cite{HW25}, the author and his collaborator obtained the exponential ergodicity in $\W_1$ with singular interaction $|h|\in \tilde{L}^k(\R^d)$ for some $k>d$, which includes the case $|h(x)|\leq c|x|^{-\beta}$ for $\beta\in(0,1)$. In \cite{HW25}, the drift $b^{(1)}$ is required to be Lipschitz continuous, which is important for the estimate of heat kernel. As far as we know, the well-posedness for \eqref{E0i} with super-linear drift $b^{(1)}$ and singular interaction $|h|\in \tilde{L}^k(\R^d)$ for some $k\in(d,\infty)$ is still open. The main difficulty is that the estimate of heat kernel is unknown when $b^{(1)}$ is of super-linear growth. Fortunately, when $k=\infty$, i.e. in the case of bounded measurable interaction, the well-posedness of \eqref{E0i} can be ensured when $b^{(1)}$ is super-linear in the sense of \eqref{mon} below, such as $b^{(1)}(x)=-|x|^2x+x$. For the exponential ergodicity in $\W_1$, we will concentrate on this bounded measurable interaction case and the regularity estimate \eqref{koz} plays an important role.

\subsection{Notations}
For $k\in(0,\infty)$, let
$$\scr P_k:=\big\{\mu\in \scr P: \mu(|\cdot|^k)<\infty\big\}.$$
Note that for $k\in (0,\infty)$, the  $L^k$-Wasserstein distance is defined as
\begin{align*}\W_k(\mu,\nu)= \inf_{\pi\in \C(\mu,\nu)}\left( \int_{\R^d\times\R^d} |x-y|^k \pi(\d x,\d y)\right)^{\frac{1}{1\vee k}},\ \ \mu,\nu\in\scr P_k.
\end{align*}
where $\C(\mu,\nu)$ is the set of all couplings of $\mu$ and $\nu$.
We will also use the total variation distance $\|\mu-\nu\|_{var}:=\sup_{|f|\leq 1}|\mu(f)-\nu(f)|$. In fact,
$$\|\mu-\nu\|_{var}=2 \inf_{\pi\in \C(\mu,\nu)} \int_{\R^d\times\R^d} 1_{\{x\neq y\}} \pi(\d x,\d y)=:2\W_0(\mu,\nu).$$

To characterize the singularity of $\sigma$ on the measure variable, we introduce the distance $\W_\eta(\eta\in(0,1])$ below instead of $\W_{\psi}$  for some $\psi$ satisfying \eqref{ESC}, which was used in \cite[Theorem 1.3]{HRWJDE}. In fact, there is no essential difference between these two distances and we use $\W_\eta$ just for simplicity.
By the dual formula, for $\eta\in(0,1]$, it holds
$$\W_\eta(\mu,\nu)= \sup_{[f]_\eta\leq 1}\big|\mu(f)-\nu(f)\big|,\ \ \mu,\nu\in\scr P_\eta,$$
where
$$[f]_\eta:=\sup_{x\neq y}\frac{|f(x)-f(y)|}{|x-y|^\eta}.$$
The SDE \eqref{E0} is called well-posed for distributions in $\scr P_k$, if for any initial value $X_0$ with $\L_{X_0}\in \scr P_k$ (respectively, any initial distribution $\gg\in \scr P_k$), it has a unique solution (respectively, a unique weak solution)
$X=(X_t)_{t\in [0,T]}$ such that $\L_{X_\cdot}:=(\L_{X_t})_{t\in [0,T]}\in C([0,T];\scr P_k).$ In this case,  for any $\gg\in \scr P_k$, let
  $P_t^*\gg=\L_{X_t^\gg}$ for the solution $X_t^\gg$ with $\L_{X_0^\gg}=\gg$. %and define
%  $$P_tf(\gamma)=\E f(X_t^\gamma)=\int_{\R^d}f(x)(P_t^\ast\gamma)(\d x), \ \ f\in\scr B_b(\R^d).$$

%The entropy-cost inequality is formulated as
%\begin{align}\label{Enp}\Ent(P_t^*\gg|P_t^*\tt\gg)\le \ff{c}t \W_2(\gg,\tt\gg)^2,\ \ t\in (0,T], \gg,\tt\gg\in \scr P_2,\end{align}
%where $\Ent$ is the relative entropy.
%Observe that \eqref{Enp} is equivalent to the log-Harnack inequality
%$$P_t \log f(\tt \gg)\le \log P_t f(\gg)+ \ff{c}{t} \W_2(\gg,\tt\gg)^2,\ \ t\in (0,T], \ f\in \B_b^+(\R^d), \ \gg,\tt\gg\in \scr P_2,$$
%one can refer to
%  \cite{HW, HW22+,HW22b, FYW1, FYW2, W10,YZ} and the monograph  \cite{Wbook} for more models on log-Harnack inequalities.
%, i.e. for any $\mu,\nu\in \scr P$,
%$\Ent(\nu|\mu):=\infty$ if $\nu$ is not absolutely continuous with respect to $\mu$, while
%$$\Ent(\nu|\mu):=\mu(\rr\log\rr)=\int_{\R^d}(\rr\log\rr)\d\mu,\ \ \text{if}\ \rr:= \ff{\d\nu}{\d\mu}\ \text{exists. }$$
%By Pinsker's inequality, we derive from \eqref{Enp} that

%with

%Although the well-posedness has been proved in \cite{HWJMAA, Ren, FYW2} under even weaker condition where $b_t^{(1)}(x,\cdot)$ is Lipschitz in the sum of  a weighted variation distance and the Wasserstein distances, regularity estimates are not yet available in such  singular situation.

\subsection{Contributions and difficulties}

The first contribution of this paper is to derive the regularity estimate \eqref{koz}.
%\begin{align*}\W_{\vv}(P_t^\ast\gamma^1,P_t^\ast\gamma^2)\leq ct^{-\frac{1}{2}+\frac{\vv}{2}}\W_{1}(\gamma^1,\gamma^2),\ \ \gamma^1,\gamma^2\in\scr P_1, t\in(0,T], \vv\in[0,1).\end{align*}
We will consider two cases respectively.

In the first case, the interaction in the drift is merely assumed to be bounded measurable while diffusion term is allowed to be Lipschitz continuous under $\W_\eta$($\eta\in(0,1]$) in the measure  variable.
When a function defined on $\scr P$ is Lipschitz continuous under $\W_{\eta}$ with $\eta\in(0,1)$ or $\|\cdot\|_{var}$, the function is actually singular in measure variable. This is quite different from the Lipschitz continuity in $\W_k$$(k\geq 1)$. As we know, when the coefficients are Lipschitz continuous under $\W_k$$(k\geq 1)$, the synchronous coupling is available. However, the synchronous coupling is invalid if the coefficients are Lipschitz continuous under $\W_{\eta}$ with $\eta\in(0,1)$ or $\|\cdot\|_{var}$. Some other techniques have to be carried out. When the diffusion term is distribution free, Girsanov's transform can be applied to estimate $\W_{\eta}$ with $\eta\in(0,1)$ or $\|\cdot\|_{var}$ between two solutions to SDEs with different drifts, see for instance \cite{HW22b,L,FYW3}.
However, the technique of the Girsanov transform is unavailable
in the case that the diffusion coefficients are distribution dependent. To overcome this essential difficulty, we adopt Duhamel's formula, basing on which some crucial estimates with respect to $\W_{\eta}$ with $\eta\in(0,1)$ or  $\|\cdot\|_{var}$ for solutions of SDEs with different diffusion coefficients are established. Compared with a series in the method of parametrix expansion of \cite{CR,CF,HWJMAA}, we make estimate directly basing on Duhamel's formula and this procedure becomes more simplified.

In the second case, the drift contains two parts: the distribution free $b^{(1)}$ satisfies monotone condition and the interacting term $b^{(0)}$ owns bounded measurable interaction. Unlike the first case, the main difficulty lies in that the Duhamel formula is unknown. We will adopt the Yosida approximation as well as approximation by mollifier to overcome this main difficulty and obtain the regularity estimate \eqref{koz}.

The second contribution is the exponential ergodicity in $\W_1$ for the second case, where some partially dissipative assumption is presented for $b^{(1)}$. Since the interaction is only bounded and measurable, it is required to use the regularity estimate \eqref{koz}.

The remaining of the paper is organized as follows:
In Section 2, we give the main results including the well-posedness, regularity estimate in Theorem \ref{EEU} and the exponential ergodicity in Theorem \ref{T0}; The proofs for Theorem \ref{EEU} and Theorem \ref{T0} will be offered in Section 3 and Section 4 respectively. In Section 5, we provide the gradient estimate as well as a comparison estimate for two classical SDEs with monotone drift plus a bounded measurable one, which is crucial in the proof of the regularity estimate \eqref{koz} in the second case and of independent interest.

Throughout this paper, $c$ or $C$ denote a constant and the values may change from one appearance to another one.
\section{Main results}
\subsection{Regularity estimate}
In this part, we fix $T>0$ and assume $b_t(x,\mu)=b_t^{(1)}(x)+b_t^{(0)}(x,\mu)$. Let
$$\|f\|_{\tilde{L}_p^q(T)}:= \sup_{z\in\R^d}\bigg\{\int_{0}^T \bigg(\int_{|x-z|\le 1}|f_u(x)|^p\d x\bigg)^{\ff q p} \d u\bigg\}^{\ff 1q}.$$
%Let
%$$\|f\|_{\tilde{L}_p^q(T)}:= \sup_{z\in\R^d}\bigg\{\int_{0}^T \bigg(\int_{|x-z|\le 1}|f_u(x)|^p\d x\bigg)^{\ff q p} \d u\bigg\}^{\ff 1q}.$$
We make the following assumptions.
\begin{enumerate}
\item[(A1)] $(\si_t\si^*_t)(x,\gamma)$ is invertible. There exists a constant $K>1$, $\beta\in(0,1]$ and $0\leq f\in \tilde{L}_p^q(T)$ with $\frac{d}{p}+\frac{2}{q}<1$ such that
 $$K^{-1}I_{d\times d}\leq \si\si^*\leq KI_{d\times d},\ \ \|\nabla\sigma_t(x,\gamma)\|\leq f_t(x),$$
     $$\|\sigma_t(x,\gamma) -\sigma_t(y,\gamma)\|\leq K|x-y|^{\beta},\ \ \gamma\in \scr P,t\in[0,T],x,y\in\R^d.$$
\item[(A2)] $b^{(0)}$ is bounded and
\beg{align*} |b^{(0)}_t(x,\gamma^1)-b^{(0)}_t(x,\gamma^2)|\leq K\|\gamma^1-\gamma^2\|_{var},\ \  t\in[0,T], x\in\R^d, \gamma^1,\gamma^2\in \scr P.\end{align*}
\item[(A3)] Let $K$, $\beta$ be in (A1). There exists a constant $\eta\in(0,1]$ such that for any $ t\in[0,T], \gamma^1,\gamma^2\in \scr P_\eta,$
$$\|\sigma_t(x,\gamma^1) -\sigma_t(x,\gamma^2)\|\leq K\W_{\eta}(\gamma^1,\gamma^2),$$
and
    \begin{align*}&\sup_{x\neq y}\frac{|[(\sigma_t\sigma_t^\ast)(x,\gamma^1)-(\sigma_t\sigma_t^\ast)(x,\gamma^2)]- [(\sigma_t\sigma_t^\ast)(y,\gamma^1)-(\sigma_t\sigma_t^\ast)(y,\gamma^2)]|}{|x-y|^{\beta}}\leq K\W_{\eta}(\gamma^1,\gamma^2).
\end{align*}
\item[(A4)] $b^{(1)}$ is locally bounded. There exists a constant  $K_1\geq 0$ such that
      $$\<b_t^{(1)}(x)-b^{(1)}_t(y),x-y\>\leq K_1|x-y|^2,\ \ t\in[0,T], x,y\in\R^d. $$
%      and $$\|\sigma(x)-\sigma(y)\|_{HS}^2\leq K_1|x-y|^2$$
\end{enumerate}
\begin{exa}\label{exabo}
Let $m=d$,  $h:\R^d\to\R^d$ be bounded and measurable and $g:\R^d\to\R^d\otimes \R^d$ be bounded and satisfy
$$|g(x)-g(y)|\leq c|x-y|^{\eta},$$
and
$$|[(gg^\ast)(x-z)-(gg^\ast)(x-z')]-[(gg^\ast)(y-z)-(gg^\ast)(y-z')]|\leq c|x-y|^\beta|z-z'|^\eta$$
for some constant $c>0$. Then $$b_t^{(0)}(x,\gamma)=\int_{\R^d}h(x-y)\gamma(\d y)$$ satisfies $(A2)$ and $$\sigma_t(x,\gamma)=\sqrt{K_0I_{d\times d}+\int_{\R^d}(gg^\ast)(x-y)\gamma(\d y)}$$satisfies $(A3)$.
%$$|[g(x-z)-g(x-z')]-[g(y-z)-g(y-z')]|\leq c|x-y|^\beta|z-z'|^\eta.$$

%$$\tilde{\sigma}(x,\mu)=\int_{\R^d}h(x-y)\mu(\d y)$$
%and set
%$$\sigma(x,\mu)=\sqrt{KI_{d\times d}+\tilde{\sigma}(x,\mu)}$$
\end{exa}
\begin{thm}\label{EEU} Assume (A1) and (A2). Then the following assertions hold.

\begin{enumerate}
\item[(i)]Suppose $b^{(1)}=0$ and (A3) holds. Then \eqref{E0} is well-posed in $\scr P_\eta$ and for any $k\geq 1$, there exists a constant $c(k)>0$ such that
    \begin{align}\label{tyq}
    (P_t^\ast\gamma)(|\cdot|^k)\leq c(k)(1+\gamma(|\cdot|^k)),\ \ t\in[0,T],\gamma\in\scr P_k.
    \end{align}
 Moreover, the following regularity estimate holds: for any $\vv\in[0,1]$,
\begin{equation}\label{REG}\begin{split}
&\W_{\vv}(P_t^\ast\gamma^1,P_t^\ast\gamma^2)\leq ct^{-\frac{1}{2}+\frac{\vv}{2}}\W_{1}(\gamma^1,\gamma^2),\ \ \gamma^1,\gamma^2\in\scr P_1, t\in(0,T].
\end{split}\end{equation}

\item[(ii)] Suppose $\sigma_t(x,\mu)=\sigma_t(x)$ and (A4) holds. Then \eqref{E0} is well-posed in $\scr P$ and \eqref{tyq}-\eqref{REG} hold with constant $c$ increasing in $\|b^{(0)}\|_{\infty}$.
\end{enumerate}
\end{thm}
 \begin{rem}\label{WEL} (1) For the well-posedness of \eqref{E0}, the diffusion coefficient $\sigma$ in \cite{CF} is assumed to have $\eta$-H\"{o}lder continuous linear functional derivative on the measure variable, which yields that $\sigma$ is Lipschitz continuous under the distance
 $$\widetilde{\W}_{\eta}(\mu,\nu):=\inf_{\pi\in\C(\mu,\nu)}\int_{\R^d\times\R^d}(|x-y|^\eta\wedge 1)\pi(\d x,\d y).$$
 In (A3), instead of introducing the linear functional derivative, $\sigma$ is allowed to be Lipschitz continuous under $\W_\eta$. Since
 $\widetilde{\W}_{\eta}\leq \W_{\eta}\wedge 1$,
 the conditions in (A3) are weaker than those in \cite{CF}.

(2) In \cite{HWJMAA}, the well-poseness is ensured when $\sigma$ is Lipschitz continuous under $\W_k$ for some $k\geq 1$ and (A3) is replaced by
 \begin{align*}&\frac{|[(\sigma_t\sigma_t^\ast)(x,\gamma^1)-(\sigma_t\sigma_t^\ast)(y,\gamma^1)]- [(\sigma_t\sigma_t^\ast)(x,\gamma^2)-(\sigma_t\sigma_t^\ast)(y,\gamma^2)]|}{|x-y|}\leq K\W_k(\gamma^1,\gamma^2).
\end{align*}
In the present Theorem \ref{EEU}(i), $\sigma$ is Lipschitz under $\W_{\eta}$ with $\eta\in(0,1)$, which allows $\sigma$ to be singular on the measure variable.

(3) Compared with \cite[Theorem 1.3(2)]{HRWJDE} or \cite[Theorem 1.2]{HWJMAA}, where the drift is Lipschitz continuous in the measure variable under $\W_\psi$(with some increasing and concave $\psi$ satisfying $\lim_{t\to0}\psi(t)=0$) or $\W_k$
(with $k\geq 1$), to derive the regularity estimate \eqref{REG} in Theorem \ref{EEU}(i), the drift $b$ is allowed to be Lipschitz continuous in the measure variable under total variation distance $\|\cdot\|_{var}$. This means that the interaction in the drift is merely required to be bounded and measurable, see Example \ref{exabo} above.
   \end{rem}
   \subsection{Exponential ergodicity in $\W_1$}

As an application of the regularity estimate \eqref{REG} in Theorem \ref{EEU}(ii), in this part, we investigate the exponential ergodicity in $\W_1$ for \eqref{E0i}. To this end, we assume $b(x,\mu)=b^{(1)} (x)+b^{(0)} (x, \mu)$ and make the following assumptions on $b^{(1)}, b^{(0)}$ and $\si$.

\beg{enumerate} \item[{\bf (C)}]  Let $a=\sigma\sigma^\ast$. There exist constants $\delta_1\geq \delta_2>0$ such that
  $$ \delta_2I_{d\times d}\leq a\leq \delta_1I_{d\times d}.$$
 Moreover, $b^{(1)}$ is locally bounded and there exist constants  $K_1,K_2,R\in (0,\infty)$ such that
\begin{align}\label{mon}\<b^{(1)}(x)-b^{(1)}(y),x-y\>\leq -K_2|x-y|^21_{\{|x-y|>R\}}+K_1|x-y|^21_{\{|x-y|\leq R\}},
\end{align}
      and $$\|\sigma(x)-\sigma(y)\|_{HS}^2\leq K_1|x-y|^2.$$
There exists a constant $\eta>0$ such that
\begin{align}\label{bbd}\|b^{(0)}\|_\infty\le \eta,\ \ \  \|b^{(0)}(\cdot,\mu)-b^{(0)}(\cdot,\nu)\|_\infty\le \eta\|\mu-\nu\|_{var},\ \ \ \mu,\nu\in \scr P.
\end{align}
 \end{enumerate}

Under {\bf(C)}, by Theorem \ref{EEU}(ii), \eqref{E0i} is well-posed and \eqref{tyq}-\eqref{REG} hold.
%\begin{equation*} \E\bigg[\sup_{t\in [0,T]} |X_t|^n\bigg]<\infty,\ \ \ \text{if}\ \E|X_0|^n<\infty,\ \ \ T\in (0,\infty),\ n\in\mathbb N.
%\end{equation*}
%Moreover,   for any   $n\geq 1, T>0$, there exists a constant $c(n,T)$ such that
%\beq\label{NES13} \E\bigg[\sup_{t\in [0,T]} |X_t|^n\bigg|\F_0\bigg]\le c(n,T) (1+|X_0|^n), \ \ t\in[0,T].\end{equation}
The following result shows that  \eqref{E0i} is exponential ergodic in $\W_1$ and $\|\cdot\|_{var}$  when $\eta$ is small enough.

\begin{thm}\label{T0} Assume {\bf (C)}.
 There  exists a constant $\eta_0>0$ depending on  $\delta_1,\delta_2, K_1,K_2,R$ such that   if  $\eta\le \eta_0$,   \eqref{E0i} has a unique invariant probability $\mu_\infty$
and
\beq\label{EX1}
 \W_1(P_t^\ast\mu,\mu_\infty)\le c \e^{-\ll t}\, \W_1(\mu,\mu_\infty),\ \ \ \mu\in \scr P_1,\ t\ge 0,\end{equation}
 \beq\label{EX2}  \|P_t^*\mu-\mu_\infty \|_{var}\le c \e^{-\ll t} (t\wedge1)^{-\ff 1 2} \W_1(\mu,\mu_\infty),\ \ \ \mu\in \scr P_1,\ t> 0\end{equation}
hold for some constants $c,\ll\in (0,\infty).$
\end{thm}
\begin{rem}\label{sinke}
In \cite[Theorem 3.1]{WSB}, the exponential ergodicity in total variation distance is also derived if
$$\<b^{(1)}(x),x\>\leq \Phi(|x|^2),\ \ x\in\R^d$$
for some $\Phi:[0,\infty)\to[1,\infty)$ satisfying $\int_0^\infty\frac{\d s}{\Phi(s)}<\infty$. In the present Theorem \ref{T0}, it is allowed that $\Phi(x)=x$.
\end{rem}
\section{Proof of Theorem \ref{EEU}}
\subsection{Proof of Theorem \ref{EEU}(i)}
\subsubsection{Duhamel's formula and some crucial estimates}
This part comes from \cite[Lemma 3.1]{HWJMAA}.
 For any $z\in\R^d,\nu\in C([0,T];\scr P_\eta)$, the weakly continuous map from $[0,T]$ to $\scr P_\eta$, let
\begin{equation*}
  a_{s,t}^{z,\nu} :=  \int_s^t(\si_u\si_u^*)(z, \nu_u)\d u,\ \ 0\leq s\leq t\leq T, \end{equation*}
%Obviously, $(A2)$   implies
%\beq\label{mv}\beg{split}
%&|m_{s,r}^{\mu,z}-m_{s,r}^{\nu,z}|\leq K(t-s) \|\mu-\nu\|_{s,r,\theta,TV},\\
% &\|a_{s,r}^{\mu,S}-a_{s,r}^{\nu,S}\|\le K\W_{s,r,\theta}\big(\Gamma(\mu), \Gamma(\nu)\big)(S_r-S_s),\\
% &\|a_{s,t}^{z,\nu^1}-a_{s,t}^{z,\nu^2}\|\le K\int_s^t\left[\sum_{i=1}^N\W_{\eta_i}\big(\nu_u^1,\nu_u^2\big)+\W_{k}\big(\nu_u^1,\nu_u^2\big)\right]\d u,\\
% &\frac{1}{K(t-s)}\leq \|[a_{s,t}^{z,\nu}]^{-1}\|\leq \frac{K}{t-s},\ \ \nu^1,\nu^2,\nu\in C([0,T];\scr P_k).\end{split}
%\end{equation}
and
\begin{equation*} q_{s,t}^{z,\nu}(x,y)=\ff{\exp[-\ff 1 {2} \<(a_{s,t}^{z,\nu})^{-1}(y-x),  y-x\>]}{(2\pi )^{\ff d 2} ({\rm det} \{a_{s,t}^{z,\nu}\})^{\ff 1 2}},\ \ x,y\in \R^d, 0\leq s< t\leq T,\end{equation*}
which is the transition density function of the Gaussian process:
\begin{align*}Y_{s,t}^{x,z,\nu}:= x + \int_s^t \si_u(z,\nu_u)\d W_u,\ \ 0\leq s\leq t\leq T, x\in\R^d.\end{align*}
Recall that $K$ is in (A1) and set
\begin{align*}\label{tiq}\tilde{q}_{s,t}(x,y)=\ff{\exp[-\frac{|y-x|^2}{4K(t-s)}]}{(4K\pi(t-s))^{\ff d 2}},\ \ x,y\in \R^d, 0\leq s< t\leq T,\end{align*}
the transition density of $\sqrt{2K}W_t$.
%It is clear that
%\begin{equation}\label{NNY}
%\E\left(\int_{\R^d}\tilde{q}_{s,t}^{S}(x,y)(1+|y|^k)\d y\right)\leq c_0(1+|x|^k+\E(S_{t-s}^{\frac{k}{2}}))\leq c(1+|x|^k).
%\end{equation}
%It is not difficult to see that
%\beq\label{G2}\beg{split} &q_{s,t}^{z,\nu}(x,y) \bigg(1+ \ff{|x-y|^4}{(t-s)^2} \bigg) \le c_1   \tilde{q}_{s,t}(x,y),\ \ x, y,z\in\R^d, 0\le s< t\le T.\end{split}\end{equation}
The following lemma is from \cite[Lemma 3.1]{HWJMAA}.
\beg{lem}\label{L20}  Assume $(A1)$. Then there exists a constant $c>0$ such that for any $0\leq s<t\le T, x,y\in \R^d$ and $
\nu^1,\nu^2,\nu\in C([0,T];\scr P_\eta), i=0,1,2$,
\beq\label{g1'}|\nabla^i q_{s,t}^{z,\nu}(\cdot,y)(x)|\leq c (t-s)^{-\frac{i}{2}} \tilde{q}_{s,t}(x,y),\end{equation}
\beq\label{g2'}
 |\nabla^i q_{s,t}^{z,\nu^1}(\cdot,y)(x)-\nabla^i q_{s,t}^{z,\nu^2}(\cdot,y)(x)|\leq   c(t-s)^{-\frac{i}{2}} \tilde{q}_{s,t} (x,y)  \frac{\int_s^t\W_{\eta}\big(\nu_u^1,\nu_u^2\big)\d u}{t-s}.
\end{equation}
\end{lem}
 % Recall the map
%\beq\label{PI} \Phi_{s_0,\cdot}^\gg: C([s_0,T];\scr P_k)\to C([s_0,T];\scr P_k);\ \
 %\Phi^\gg_{s_0,t}(\mu):=\L_{X_{s_0,t}^{\gamma,\mu}},\ \   t\in [s_0,T], \mu\in C([s_0,T];\scr P_k).\end{equation}
%To prove the existence and uniqueness of the fixed point for $\Phi_{s_0,\cdot}^\gg$, we  investigate  the contraction of this map with respect to the complete metric
%$$\|\mu-\nu\|_{var,k,\delta}:= \sup_{t\in [0,T]}\e^{-\delta t} \{\|\mu_t-\nu_t\|_{k,var}+\W_k(\mu_t,\nu_t)\},\ \ \mu,\nu\in C([s_0,T];\scr P_k).$$
\subsubsection{Duhamel's formula}
For any $\nu\in C([0,T];\scr P_\eta)$, consider
 \beq\label{Et}\beg{split} & \d  {X}_{s,t}^{x,\nu}= b_t( {X}_{s,t}^{x, \nu}, \nu_t)\d t+\sigma_t({X}_{s,t}^{x, \nu},\nu_t)\d W_{t},\ \ 0\leq s\leq t\leq T,\ \ X_{s,s}^{x,\nu}=x\in\R^d.\end{split}\end{equation}
Then under (A1)-(A2), \eqref{Et} has  a unique solution and let $p_{s,t}^{\nu}(x,\cdot),0\leq s<t\leq T$ be the distribution density function of ${X}_{s,t}^{x,\nu}$.
%Let
%\begin{align}\label{LLT}&P_{s,t}^{\mu,\nu}f(x)=\E f(X_{s,t}^{x,\mu,\nu})=\int_{\R^d}p_{s,t}^{\mu,\nu}(x,y)f(y)\d y,\nonumber\\
%&Q_{s,t}^{\nu} f(x)=\int_{\R^d}q_{s,t}^{z,\nu}(x,z)f(z)\d z,\ \ f\in \scr B_b(\R^d).
% \end{align}
Let ${X}_{s,t}^{\gamma,\nu}$ be the solution to \eqref{Et} from initial distribution $\gamma\in\scr P$. For simplicity, we denote ${X}_{t}^{\gamma,\nu}={X}_{0,t}^{\gamma,\nu}$ and $a_t=\sigma_t\sigma_t^\ast$. Under (A1)-(A2), the Duhamel formula holds, i.e.
\begin{align*}
    &p_{s,t}^{\nu}(x,z)\nonumber\\
    &= q_{s,t}^{z,\nu}(x,z)+\int_s^t \int_{\R^d}p_{s,r}^{\nu}(x,y) \left\<b_r(y, \nu_r),\nabla q_{r,t}^{z,\nu}(\cdot,z)(y)\right\>\d y\d r\\
    &+\frac{1}{2}\int_s^t  \int_{\R^d}p_{s,r}^{\nu}(x,y) \mathrm{tr}\{[a_r(y,\nu_r)-a_r(z,\nu_r)]\nabla^2 q_{r,t}^{z,\nu}(\cdot,z)(y)\}\d y\d r\nonumber,
\end{align*}
which yields for any $f\in \scr B_b(\R^d)$, $\gamma\in\scr P$,
\begin{align}\label{PTf}
    \nonumber&\int_{\R^d}\int_{\R^d} p_{s,t}^{\nu}(x,z)f(z)\d z\gamma(\d x)\\
    &=\int_{\R^d}\int_{\R^d} q_{s,t}^{z,\nu}(x,z)f(z)\d z\gamma(\d x)+\int_{\R^d}R_{s,t}^{\nu}f(x)\gamma(\d x)
    \end{align}
    with
\begin{align*}
    \nonumber&R_{s,t}^{\nu}f(x):=\int_s^t \int_{\R^d}p_{s,r}^{\nu}(x,y) \left\<b_r(y, \nu_r),\int_{\R^d}\nabla q_{r,t}^{z,\nu}(\cdot,z)(y)f(z)\d z\right\>\d y\d r\\
    &+\frac{1}{2}\int_s^t  \int_{\R^d}p_{s,r}^{\nu}(x,y)\mathrm{tr}\left\{\int_{\R^d}[a_r(y,\nu_r)-a_r(z,\nu_r)]\nabla^2 q_{r,t}^{z,\nu}(\cdot,z)(y)f(z)\d z\right\}\d y\d r.
\end{align*}
%\begin{align}\label{PTf}
%    &\int_{\R^d}P_{s,t}^{\mu,\nu}f(x)\gamma(\d x)\nonumber\\
%    & =\int_{\R^d}Q_{s,t}^{\nu} f(x)\gamma(\d x)+\int_{\R^d}\int_s^t  P_{s,r}^{\mu,\nu} \left\<b_r(\cdot, \mu_r),\nabla Q_{r,t}^{\nu}f\right\>(x)\d r\gamma(\d x)\\ \nonumber
%    &+\frac{1}{2}\int_{\R^d}\int_s^t  P_{s,r}^{\mu,\nu}\int_{\R^d} \mathrm{tr}\{[(\sigma_r\sigma_r^\ast)(\cdot,\nu_r)-(\sigma_r\sigma_r^\ast)(z,\nu_r)]\nabla^2 q_{r,t}^{z,\nu}(\cdot,z)f(z)\d z\}(x)\d r\gamma(\d x)\\
%    &\nonumber=:\int_{\R^d}Q_{s,t}^{\nu} f(x)\gamma(\d x)+\int_{\R^d}R_{s,t}^{\mu,\nu}f(x)\gamma(\d x).
%\end{align}

\subsubsection{Estimate of $\|\L_{{X}_{s,t}^{\gamma,\nu^1}}-\L_{{X}_{s,t}^{\gamma,\nu^2}}\|_{var}$ }
First, we provide estimate of $\|\L_{{X}_{s,t}^{\gamma,\nu^1}}-\L_{{X}_{s,t}^{\gamma,\nu^2}}\|_{var}$ by Duhamel's formula \eqref{PTf}.
\begin{lem}\label{RST}
Assume (A1)-(A3) and $b^{(1)}=0$. Then there exists a constant $c>0$ such that
\begin{equation}\begin{split}\label{LLY}
&\|\L_{{X}_{s,t}^{\gamma,\nu^1}}-\L_{{X}_{s,t}^{\gamma,\nu^2}}\|_{var}\\
&\leq c\frac{\int_s^t\W_{\eta}(\nu^1_u,\nu^2_u)\d u}{t-s}+c\int_{s}^t\|\L_{{X}_{s,r}^{\gamma,\nu^1}}-\L_{{X}_{s,r}^{\gamma,\nu^2}}\|_{var}  (t-r)^{-1+\frac{\beta}{2}}\d r\\
&+c\int_{s}^t\|\nu^1_r-\nu^2_r\|_{var} (t-r)^{-\frac{1}{2}}\d r+c\int_s^t\W_{\eta}(\nu^1_r,\nu^2_r) (t-r)^{-1+\frac{\beta}{2}}\d r\\
&+c\int_{s}^t(t-r)^{-2+\frac{\beta}{2}}\int_r^t\W_{\eta}(\nu^1_u,\nu^2_u)\d u \d r.
\end{split}\end{equation}
Consequently, there exists constants $c>0,\lambda_0>1$ such that for any $\lambda\geq \lambda_0$,
\begin{equation}\begin{split}\label{MGY}
&\sup_{t\in[0,T]}\e^{-\lambda t}t^{\frac{1}{2}}\|\L_{{X}_{0,t}^{\gamma,\nu^1}}-\L_{{X}_{0,t}^{\gamma,\nu^2}}\|_{var}\\
&\leq c\lambda ^{-\frac{\beta\wedge \eta}{2}}\sup_{t\in[0,T]}\e^{-\lambda t}t^{\frac{1-\eta}{2}}\W_{\eta}(\nu^1_t,\nu^2_t) +c\lambda^{-\frac{1}{2}}\sup_{t\in[0,T]}\e^{-\lambda t}t^{\frac{1}{2}}\|\nu^1_t-\nu^2_t\|_{var}.
\end{split}\end{equation}
\end{lem}
\begin{proof}
Firstly, it follows from \eqref{g2'} that
\begin{align*}
&\sup_{|f|\leq 1}\left|\int_{\R^d}\int_{\R^d}[q_{s,t}^{z,\nu^1}(x,z)-q_{s,t}^{z,\nu^2}(x,z)]f(z)\d z\gamma(\d x)\right|\leq c_1\frac{\int_s^t\left(\W_{\eta}(\nu^1_u,\nu^2_u)\right)\d u}{t-s}.
\end{align*}
By \eqref{PTf}, we conclude that
\begin{align*}
&\sup_{|f|\leq 1}\left|\int_{\R^d}[R_{s,t}^{\nu^1}f(x)-R_{s,t}^{\nu^2}f(x)]\gamma(\d x)\right|\leq\sum_{i=1}^6\int_{s}^t |J_i|\d r ,
\end{align*}
where
\begin{align*}
 &J_1=\sup_{|f|\leq 1}\bigg|\int_{\R^d}\gamma(\d x)\int_{\R^d} [p_{s,r}^{\nu^1} -p_{s,r}^{\nu^2} ](x,y)\d y\\
 &\qquad\qquad\quad\times\left\<b_r(y, \nu^1_r),\int_{\R^d}\nabla q_{r,t}^{z,\nu^1}(\cdot,z)(y)f(z)\d z\right\>\bigg|,
\end{align*}
\begin{align*}
  &J_2=\frac{1}{2}\sup_{|f|\leq 1}\bigg|\int_{\R^d}\gamma(\d x)\int_{\R^d} [p_{s,r}^{\nu^1} -p_{s,r}^{\nu^2} ](x,y)\d y\\
  &\qquad\qquad\quad\times \mathrm{tr}\left\{\int_{\R^d}[a_r(y,\nu^1_r)-a_r(z,\nu^1_r)]\nabla^2 q_{r,t}^{z,\nu^1}(\cdot,z)(y)f(z)\d z\right\}\bigg|,
\end{align*}
\begin{align*}
  &J_3=\int_{\R^d}\gamma(\d x)\int_{\R^d} p_{s,r}^{\nu^2} (x,y) |b_r(y, \nu^1_r)-b_r(y, \nu^2_r)|\d y\int_{\R^d}|\nabla q_{r,t}^{z,\nu^1}(\cdot,z)(y)|\d z,
 \end{align*}
 \begin{align*}
  &J_4=\frac{1}{2}\int_{\R^d}\gamma(\d x)\int_{\R^d} p_{s,r}^{\nu^2} (x,y)\d y \\
  &\quad\times \mathrm{tr}\left\{\int_{\R^d}|[a_r(y,\nu_r^1)-a_r(z,\nu^1_r)] -[a_r(y,\nu_r^2)-a_r(z,\nu_r^2)]||\nabla^2q_{r,t}^{z,\nu^1}(\cdot,z)(y)| \d z\right\},
 \end{align*}
  \begin{align*}
  &J_5=\int_{\R^d}\gamma(\d x)\int_{\R^d} p_{s,r}^{\nu^2} (x,y)\|b\|_\infty\d y\int_{\R^d}|\nabla q_{r,t}^{z,\nu^1}(\cdot,z)(y)-\nabla q_{r,t}^{z,\nu^2}(\cdot,z)(y)| \d z,
   \end{align*}
 \begin{align*}
  &J_6=\frac{1}{2} \int_{\R^d}\gamma(\d x)\int_{\R^d} p_{s,r}^{\nu^2} (x,y)\d y\\
  &\quad\times\mathrm{tr}\left\{ \int_{\R^d}[a_r(y,\nu_r^2)-a_r(z,\nu_r^2)]|\nabla^2 q_{r,t}^{z,\nu^1}(\cdot,z)(y)-\nabla^2 q_{r,t}^{z,\nu^2}(\cdot,z)(y)| \d z\right\}.
\end{align*}
We derive from (A1)-(A2) and Lemma \ref{L20} that for any $f$ with $|f|\leq 1$,
\begin{align*}\left|\left\<b_r(y, \nu^1_r),\int_{\R^d}\nabla q_{r,t}^{z,\nu^1}(\cdot,z)(y)f(z)\d z\right\>\right|\leq c_1(t-r)^{-\frac{1}{2}},
\end{align*}
and
\begin{align*}&\left|\mathrm{tr}\left\{\int_{\R^d}[a_r(y,\nu^1_r)-a_r(z,\nu^1_r)]\nabla^2 q_{r,t}^{z,\nu^1}(\cdot,z)(y)f(z)\d z\right\}\right|\leq c_2(t-r)^{-1+\frac{\beta}{2}}.
\end{align*}
So, by the definition of $\|\cdot\|_{var}$, we conclude
\begin{align*}
J_1&\leq c_1\|\L_{{X}_{s,r}^{\gamma,\nu^1}}-\L_{{X}_{s,r}^{\gamma,\nu^2}}\|_{var}  (t-r)^{-\frac{1}{2}},
\end{align*}
and
\begin{align*}
J_2&\leq c_2\|\L_{{X}_{s,r}^{\gamma,\nu^1}}-\L_{{X}_{s,r}^{\gamma,\nu^2}}\|_{var}  (t-r)^{-1+\frac{\beta}{2}}.
\end{align*}
By (A1), (A3) and \eqref{g1'},  we arrive at
\begin{align*}
J_3&\leq c_3\|\nu^1_r-\nu^2_r\|_{var} (t-r)^{-\frac{1}{2}},
\end{align*}
and
\begin{align*}
J_4&\leq c_4\W_{\eta}(\nu^1_r,\nu^2_r) (t-r)^{-1+\frac{\beta}{2}}.
\end{align*}
%Note that $\W_{r,t,k}(\Phi_{s_0,\cdot}^\gg(\mu), \Phi_{s_0,\cdot}^\gg(\nu))\leq \W_{s_0,t,k}(\Phi_{s_0,\cdot}^\gg(\mu), \Phi_{s_0,\cdot}^\gg(\nu)), r\in[s_0,t]$.
It follows from (A1)-(A2) and \eqref{g2'} that
\begin{align*}
J_5&\leq c_5\int_r^t\W_{\eta}(\nu^1_u,\nu^2_u)\d u(t-r)^{-\frac{3}{2}},
\end{align*}
and
\begin{align*}
J_6&\leq c_6\int_r^t\W_{\eta}(\nu^1_u,\nu^2_u)\d u(t-r)^{-2+\frac{\beta}{2}}.
\end{align*}
Combining all the estimates, we derive \eqref{LLY}.

(2) Firstly, it follows from the FKG inequality that
\begin{align}\label{kpt13}
\nonumber&\sup_{t\in[0,T]}\e^{-\lambda t}t^{\frac{1}{2}}\frac{\int_0^t\W_{\eta}(\nu^1_u,\nu^2_u)\d u}{t}\\
\nonumber&\leq \sup_{u\in[0,T]}\e^{-\lambda u}u^{\frac{1-\eta}{2}}\W_{\eta}(\nu^1_u,\nu^2_u)\sup_{t\in[0,T]}t^{-\frac{1}{2}}\int_0^tu^{\frac{\eta-1}{2}}\e^{-\lambda(t-u)}\d u\\
&\leq \sup_{u\in[0,T]}\e^{-\lambda u}u^{\frac{1-\eta}{2}}\W_{\eta}(\nu^1_u,\nu^2_u)\sup_{t\in[0,T]}t^{\frac{1}{2}} \frac{1}{t}\int_0^tu^{\frac{\eta-1}{2}}\d u\frac{1}{t}\int_0^t\e^{-\lambda(t-u)}\d u\\
\nonumber&\leq\frac{2}{\eta+1} \sup_{u\in[0,T]}\e^{-\lambda u}u^{\frac{1-\eta}{2}}\W_{\eta}(\nu^1_u,\nu^2_u)\sup_{t\in[0,T]}\left[t^{-1+\frac{\eta}{2}}\frac{1-\e^{-\lambda t}}{\lambda}\right]\\
\nonumber&\leq \frac{2}{\eta+1}\lambda^{-\frac{\eta}{2}}\sup_{r\geq 0} [r^{-1+\frac{\eta}{2}}(1-\e^{-r})]\sup_{u\in[0,T]}\e^{-\lambda u}u^{\frac{1-\eta}{2}}\W_{\eta}(\nu^1_u,\nu^2_u),
\end{align}
where in the last step, we use the transformation of variable $\lambda t=r$.

Secondly, it holds
\begin{align*}
&\sup_{t\in[0,T]}\e^{-\lambda t}t^{\frac{1}{2}}\int_{0}^t\|\L_{{X}_{0,r}^{\gamma,\nu^1}}-\L_{{X}_{0,r}^{\gamma,\nu^2}}\|_{var}  (t-r)^{-1+\frac{\beta}{2}}\d r\\
&\leq \sup_{r\in[0,T]}\e^{-\lambda r}r^{\frac{1}{2}}\|\L_{{X}_{0,r}^{\gamma,\nu^1}}-\L_{{X}_{0,r}^{\gamma,\nu^2}}\|_{var}\sup_{t\in[0,T]}t^{\frac{1}{2}}\int_{0}^t  \e^{-\lambda (t-r)}r^{-\frac{1}{2}}(t-r)^{-1+\frac{\beta}{2}}\d r.
\end{align*}
By the FKG inequality, we have
\begin{align}\label{FKG13}
\nonumber&t^{\frac{1}{2}}\int_{0}^t  \e^{-\lambda (t-r)}r^{-\frac{1}{2}}(t-r)^{-1+\frac{\beta}{2}}\d r\\
&\leq t^{\frac{1}{2}}t\frac{1}{t}\int_{0}^t  r^{-\frac{1}{2}}\d r\frac{1}{t}\int_{0}^t  \e^{-\lambda (t-r)}(t-r)^{-1+\frac{\beta}{2}}\d r\\
\nonumber&\leq 2\int_{0}^t  \e^{-\lambda r}r^{-1+\frac{\beta}{2}}\d r\leq 2\lambda^{-\frac{\beta}{2}}\int_{0}^\infty  \e^{-t}t^{-1+\frac{\beta}{2}}\d t.
\end{align}
So, we get
\begin{align}\label{KKT}
\nonumber&\sup_{t\in[0,T]}\e^{-\lambda t}t^{\frac{1}{2}}\int_{0}^t\|\L_{{X}_{0,r}^{\gamma,\nu^1}}-\L_{{X}_{0,r}^{\gamma,\nu^2}}\|_{var}  (t-r)^{-1+\frac{\beta}{2}}\d r\\
&\leq 2\Gamma(\frac{\beta}{2})\lambda^{-\frac{\beta}{2}}\sup_{r\in[0,T]} \e^{-\lambda r} r^{\frac{1}{2}}\|\L_{{X}_{0,r}^{\gamma,\nu^1}}-\L_{{X}_{0,r}^{\gamma,\nu^2}}\|_{var}.
\end{align}
Similarly, we have
\begin{equation}\begin{split}\label{MMG11}
&\sup_{t\in[0,T]}\e^{-\lambda t}t^{\frac{1}{2}}\int_{0}^t\|\nu^1_r-\nu^2_r\|_{var} (t-r)^{-\frac{1}{2}}\d r\leq 2\Gamma(\frac{1}{2})\lambda^{-\frac{1}{2}}\sup_{r\in[0,T]}\e^{-\lambda r}r^{\frac{1}{2}}\|\nu^1_r-\nu^2_r\|_{var},
\end{split}\end{equation}
and
\begin{equation}\begin{split}\label{MMG110}
&\sup_{t\in[0,T]}\e^{-\lambda t}t^{\frac{1}{2}}\int_0^t\W_{\eta}(\nu^1_r,\nu^2_r) (t-r)^{-1+\frac{\beta}{2}}\d r\\
&\leq \sup_{t\in[0,T]}\e^{-\lambda t}t^{\frac{1-\eta}{2}}\W_{\eta}(\nu^1_t,\nu^2_t)\sup_{t\in[0,T]}t^{\frac{1}{2}}\int_0^t r^{-\frac{1-\eta}{2}}\e^{-\lambda(t-r) }(t-r)^{-1+\frac{\beta}{2}}\d s\\
&\leq \frac{2}{1+\eta}\sup_{t\in[0,T]}\e^{-\lambda t}t^{\frac{1-\eta}{2}}\W_{\eta}(\nu^1_t,\nu^2_t)\sup_{t\in[0,T]}t^{\frac{\eta}{2}}\int_0^t \e^{-\lambda s}s^{-1+\frac{\beta}{2}}\d s\\
&\leq \frac{2}{1+\eta}\Gamma(\frac{\beta}{2})\lambda^{-\frac{\beta}{2}}T^{\frac{\eta}{2}}\sup_{t\in[0,T]}\e^{-\lambda t}t^{\frac{1-\eta}{2}}\W_{\eta}(\nu^1_t,\nu^2_t),
\end{split}\end{equation}
where in the last step, we use the transformation of variable $\lambda s=r$.

Finally, we have
\begin{equation*}\begin{split}
&\sup_{t\in[0,T]}\e^{-\lambda t}t^{\frac{1}{2}}\int_{0}^t\int_r^t\W_{\eta}(\nu^1_u,\nu^2_u)\d u(t-r)^{-2+\frac{\beta}{2}}\d r\\
&\leq\sup_{u\in[0,T]}\e^{-\lambda u}u^{\frac{1-\eta}{2}}\W_{\eta}(\nu^1_u,\nu^2_u) \sup_{t\in[0,T]}t^{\frac{1}{2}}\int_{0}^t\int_r^t\e^{-\lambda (t-u)}u^{\frac{-1+\eta}{2}}\d u(t-r)^{-2+\frac{\beta}{2}}\d r.
\end{split}\end{equation*}
%It follows from the FKG inequality that
%\begin{align*}
%&\sup_{t\in[0,T]}t^{\frac{1}{2}}\int_{0}^t\int_r^tu^{\frac{-1+\eta}{2}}\e^{-\lambda(t-u)}\d u(t-r)^{-2+\frac{\beta}{2}}\d r\\
%&\leq \sup_{t\in[0,T]}t^{\frac{1}{2}}\int_{0}^t(t-r)\frac{1}{t-r}\int_r^tu^{\frac{-1+\eta}{2}}\d u\frac{1}{t-r}\int_r^t\e^{-\lambda(t-u)}\d u(t-r)^{-2+\frac{\beta}{2}}\d r\\
%&\leq \frac{2}{1+\eta} \sup_{t\in[0,T]}t^{\frac{1}{2}}\int_{0}^t(t^{\frac{1+\eta}{2}} -r^{\frac{1+\eta}{2}}) \frac{1-\e^{-\lambda(t-r)}}{\lambda(t-r)}(t-r)^{-2+\frac{\beta}{2}}\d r\\
%&\leq \frac{2}{1+\eta}\sup_{t\in[0,T]}t^{\frac{1}{2}}\int_{0}^ts^{\frac{-1+\eta}{2}}\frac{1-\e^{-\lambda s}}{\lambda}s^{-2+\frac{\beta}{2}}\d s
%\end{align*}
It follows from the the Fubini theorem and the FKG inequality that
\begin{align*}
&\sup_{t\in[0,T]}t^{\frac{1}{2}}\int_{0}^t\int_r^tu^{\frac{-1+\eta}{2}}\e^{-\lambda(t-u)}\d u(t-r)^{-2+\frac{\beta}{2}}\d r\\
&=\sup_{t\in[0,T]}t^{\frac{1}{2}}\int_{0}^t\int_0^u(t-r)^{-2+\frac{\beta}{2}}\d ru^{\frac{-1+\eta}{2}}\e^{-\lambda(t-u)}\d u\\
&\leq (1-\frac{\beta}{2})^{-1}\sup_{t\in[0,T]}t^{\frac{1}{2}}\int_{0}^t(t-u)^{-1+\frac{\beta}{2}}u^{\frac{-1+\eta}{2}}\e^{-\lambda(t-u)}\d u\\
&\leq(1-\frac{\beta}{2})^{-1}\sup_{t\in[0,T]}t^{\frac{1}{2}}t\frac{1}{t} \int_{0}^t(t-u)^{-1+\frac{\beta}{2}}\e^{-\lambda(t-u)}\d u\frac{1}{t}\int_{0}^tu^{\frac{-1+\eta}{2}}\d u\\
&\leq(\frac{1+\eta}{2})^{-1}(1-\frac{\beta}{2})^{-1}T^{\frac{\eta}{2}}\sup_{t\in[0,T]} \int_{0}^ts^{-1+\frac{\beta}{2}}\e^{-\lambda s}\d s\\
&\leq (\frac{1+\eta}{2})^{-1}(1-\frac{\beta}{2})^{-1}\Gamma(\frac{\beta}{2})\lambda^{-\frac{\beta}{2}}T^{\frac{\eta}{2}}.
\end{align*}
So, we conclude that
\begin{equation}\begin{split}\label{MMG221}
&\sup_{t\in[0,T]}\e^{-\lambda t}t^{\frac{1}{2}}\int_{0}^t\int_r^t\W_{\eta}(\nu^1_u,\nu^2_u)\d u(t-r)^{-2+\frac{\beta}{2}}\d r\\
&\leq(\frac{1+\eta}{2})^{-1}(1-\frac{\beta}{2})^{-1}\Gamma(\frac{\beta}{2})\lambda^{-\frac{\beta}{2}}T^{\frac{\eta}{2}}\sup_{u\in[0,T]}\e^{-\lambda u}u^{\frac{1-\eta}{2}}\W_{\eta}(\nu^1_u,\nu^2_u).
\end{split}\end{equation}
Taking $\lambda_0=1\vee(4c\Gamma(\frac{\beta}{2}))^{\frac{2}{\beta}}$, we get $2c\Gamma(\frac{\beta}{2})\lambda^{-\frac{\beta}{2}}\leq \frac{1}{2}, \lambda\geq \lambda_0$, which combined with \eqref{LLY}, \eqref{kpt13}-\eqref{MMG221} gives \eqref{MGY}.
\end{proof}
\subsubsection{Estimate of $\W_\vv(\L_{{X}_{0,t}^{\gamma,\nu^1}},\L_{{X}_{0,t}^{\gamma,\nu^2}})$ with $\vv\in(0,1)$}
\begin{thm}\label{WET} Assume (A1)-(A3) and $b^{(1)}=0$ . Then there exists a constant $c(\eta)\in(0,\infty)$ such that for any $\lambda\geq \lambda_0$ and $\vv\in(0,1)$,
\begin{align}\label{kpt}\nonumber&\sup_{t\in[0,T]}\e^{-\lambda t}t^{\frac{1-\vv}{2}}\W_\varepsilon(\L_{{X}_{0,t}^{\gamma,\nu^1}},\L_{{X}_{0,t}^{\gamma,\nu^2}})\\
&\leq c(\eta)\lambda ^{-\frac{\beta\wedge \eta}{2}}\sup_{t\in[0,T]}\e^{-\lambda t}t^{\frac{1-\eta}{2}}\W_{\eta}(\nu^1_t,\nu^2_t) +c(\eta)\lambda^{-\frac{1}{2}}\sup_{t\in[0,T]}\e^{-\lambda t}t^{\frac{1}{2}}\|\nu^1_t-\nu^2_t\|_{var}.
\end{align}
\end{thm}
\begin{proof}
By \cite[(1.6)]{HRWJDE}, there exists a constant $c>0$ such that
\begin{align*}
\nonumber \W_1(\L_{{X}_{0,t}^{\gamma,\nu^1}},\L_{{X}_{0,t}^{\gamma,\nu^2}})
     & \le  c\left(\int_0^t \W_{\eta}(\nu^1_r,\nu_r^2)^{2}\d r\right)^{\frac{1}{2}}+c\int_0^t\|\nu^1_r-\nu_r^2\|_{var}\d r,\ \ 0\leq t\leq T.
\end{align*}
So, we have
\begin{align*}&\sup_{t\in[0,T]}\e^{-\lambda t}\W_1(\L_{{X}_{0,t}^{\gamma,\nu^1}},\L_{{X}_{0,t}^{\gamma,\nu^2}})\\
     & \le    c\sup_{t\in[0,T]}\e^{-\lambda t}\left(\int_0^t \W_{\eta}(\nu^1_r,\nu_r^2)^{2}\d r\right)^{\frac{1}{2}}+c\sup_{t\in[0,T]}\e^{-\lambda t}\int_0^t\|\nu^1_r-\nu_r^2\|_{var}\d r\\
     & \le    c\sup_{t\in[0,T]}\e^{-\lambda t}t^{\frac{1-\eta}{2}}\W_{\eta}(\nu^1_t,\nu_t^2)\sup_{t\in[0,T]}\left(\int_0^t \e^{-2\lambda (t-u)} u^{-1+\eta}\d u\right)^{\frac{1}{2}}\\
     &+c\sup_{r\in[0,T]}\e^{-\lambda r}r^{\frac{1}{2}}\|\nu^1_r-\nu_r^2\|_{var}\sup_{t\in[0,T]}\int_0^t\e^{-\lambda (t-r)}r^{-\frac{1}{2}}\d r.
\end{align*}
By the FKG inequality, we have
\begin{align*}
&\int_0^t\e^{-\lambda (t-r)}r^{-\frac{1}{2}}\d r\leq t\frac{1}{t}\int_0^t\e^{-\lambda (t-r)}\d r\frac{1}{t}\int_0^tr^{-\frac{1}{2}}\d r\leq 2\lambda^{-\frac{1}{2}}\sup_{y\in[0,\infty)}y^{-\frac{1}{2}}(1-\e^{-y}),
\end{align*}
and
\begin{align*}
\int_0^t \e^{-2\lambda (t-u)} u^{-1+\eta}\d u
&\leq t\frac{1}{t}\int_0^t\e^{-2\lambda (t-u)}\d u\frac{1}{t}\int_0^tu^{-1+\eta}\d r\leq \frac{\eta}{1}\lambda^{-\eta}\sup_{y\in[0,\infty)}y^{\eta-1}\frac{(1-\e^{-2y})}{2}.
\end{align*}
This implies that
\begin{align}\label{W13}\nonumber&\sup_{t\in[0,T]}\e^{-\lambda t}\W_1(\L_{{X}_{0,t}^{\gamma,\nu^1}},\L_{{X}_{0,t}^{\gamma,\nu^2}})\\
&\leq c_1\lambda^{-\frac{\eta}{2}}\sup_{t\in[0,T]}\e^{-\lambda t}t^{\frac{1-\eta}{2}}\W_{\eta}(\nu^1_t,\nu_t^2)+c_1\lambda^{-\frac{1}{2}}\sup_{r\in[0,T]}\e^{-\lambda r}r^{\frac{1}{2}}\|\nu^1_r-\nu_r^2\|_{var}.
\end{align}
By \cite[Lemma 2.1]{HRWJDE}, for any $\vv\in(0,1]$, it holds
\begin{align}\label{kpz}t^{\frac{1-\vv}{2}}\W_{\vv}(\gamma,\tilde{\gamma})\leq \ss d t^{\frac{1}{2}}\|\gamma-\tilde{\gamma}\|_{var}+d\W_{1}(\gamma,\tilde{\gamma}),\ \ \gamma,\tilde{\gamma}\in\scr P_1.
\end{align}
So, we derive \eqref{kpt} from \eqref{W13}, \eqref{MGY} and \eqref{kpz}.
\end{proof}

%Since $b,\sigma$ are bounded due to (A1)-(A2), for any solution $X_t$ to \eqref{EMU}, we have
%\begin{align}\label{CTY}\E\sup_{t\in[0,T]}|X_t|^k\leq c_0(\gamma(|\cdot|^k)+T^k+T^{\frac{k}{2}})\leq c_1\gamma(1+|\cdot|^k) \end{align}
%for some constant $c_1\geq 1.$
\subsubsection{Proof of \eqref{REG}}
\begin{proof}[Proof of \eqref{REG}]
(1) Let $\nu^i \in C([0,T];\scr P_\eta)$, $\gg\in \scr P_\eta, i=1,2$.
Combining \eqref{MGY} with Theorem \ref{WET}, we can find some constant $c(\eta)>0$ such that
\begin{align}\label{cty13}\nonumber&\sup_{t\in[0,T]}\e^{-\lambda t}\left(t^{\frac{1-\vv}{2}}\W_{\vv}(\L_{{X}_{0,t}^{\gamma,\nu^1}}, \L_{{X}_{0,t}^{\gamma,\nu^2}})+t^{\frac{1}{2}}\|\L_{{X}_{0,t}^{\gamma,\nu^1}}- \L_{{X}_{0,t}^{\gamma,\nu^2}}\|_{var}\right)\\
&\leq c(\eta)\lambda ^{-\frac{\beta\wedge \eta}{2}}\sup_{t\in[0,T]}\e^{-\lambda t}(t^{\frac{1-\eta}{2}}\W_{\eta}(\nu^1_t,\nu^2_t) +t^{\frac{1}{2}}\|\nu^1_t-\nu^2_t\|_{var}),\ \ \lambda\geq\lambda_0,\vv\in(0,1).
\end{align}
Taking $\vv=\eta$, we derive
\begin{align}\label{cty}\nonumber&\sup_{t\in[0,T]}\e^{-\lambda t}\left(t^{\frac{1-\eta}{2}}\W_{\eta}(\L_{{X}_{0,t}^{\gamma,\nu^1}}, \L_{{X}_{0,t}^{\gamma,\nu^2}})+t^{\frac{1}{2}}\|\L_{{X}_{0,t}^{\gamma,\nu^1}}- \L_{{X}_{0,t}^{\gamma,\nu^2}}\|_{var}\right)\\
&\leq c(\eta)\lambda ^{-\frac{\beta\wedge \eta}{2}}\sup_{t\in[0,T]}\e^{-\lambda t}(t^{\frac{1-\eta}{2}}\W_{\eta}(\nu^1_t,\nu^2_t) +t^{\frac{1}{2}}\|\nu^1_t-\nu^2_t\|_{var}),\ \ \lambda\geq\lambda_0.
\end{align}
Taking $\lambda_1=\lambda_0\vee(2c(\eta))^{\frac{2}{\beta\wedge\eta}}$ in \eqref{cty}, we have
\begin{align*}\nonumber&\sup_{t\in[0,T]}\e^{-\lambda t}\left(t^{\frac{1-\eta}{2}}\W_{\eta}(\L_{{X}_{0,t}^{\gamma,\nu^1}}, \L_{{X}_{0,t}^{\gamma,\nu^2}})+t^{\frac{1}{2}}\|\L_{{X}_{0,t}^{\gamma,\nu^1}}- \L_{{X}_{0,t}^{\gamma,\nu^2}}\|_{var}\right)\\
&\leq \frac{1}{2}\sup_{t\in[0,T]}\e^{-\lambda t}(t^{\frac{1-\eta}{2}}\W_{\eta}(\nu^1_t,\nu^2_t) +t^{\frac{1}{2}}\|\nu^1_t-\nu^2_t\|_{var}),\ \ \lambda\geq \lambda_1.
\end{align*}
This implies that the map $\Psi^{\gamma}$ defined by $\Psi^{\gamma}_t(\nu)=\L_{{X}_{0,t}^{\gamma,\nu}},t\in[0,T]$ is strictly contractive on complete metric space $(C([0,T],\scr P_\eta),\rho_\lambda)$ with $$\rho_\lambda(\nu^1,\nu^2):=\sup_{t\in[0,T]}\e^{-\lambda t}(t^{\frac{1-\eta}{2}}\W_{\eta}(\nu^1_t,\nu^2_t) +t^{\frac{1}{2}}\|\nu^1_t-\nu^2_t\|_{var}),\ \ \nu^1,\nu^2\in C([0,T],\scr P_\eta), \lambda\geq \lambda_1.$$
By the Banach fixed point theorem, we conclude that \eqref{E0} is well-posed in $\scr P_\eta$.

(2) For any  $\gamma^i \in\scr P_{\eta}$, $i=1,2,$ let
    $$\mu_t^i=P_t^\ast\gamma^i, \ \ i=1,2, t\in [0,T].$$
    Then it holds \begin{align*}\mu_t^i=\L_{X_{0,t}^{\gamma^i,\mu^i}}, \ \ i=1,2, t\in [0,T].
    \end{align*}
%We have
%\beq\label{RS} \mu_t^i= \int_{\R^d} \L_{X_{0,t}^{x,\mu^i,\mu^i}}\gamma^i(\d x),\ \ t\in [0,T], i=1,2.\end{equation}
% \beq\label{RS*} \Phi_t^{\mu_0^2}\mu^2= P_t^*\mu_0^2,\ \ \ \textcolor{red}{\Phi_t^{\mu_0^2}\mu^1= \int_{\R^d} \L_{X_t^{x,\mu^1}}\mu_0^2(\d x)}.\end{equation}
So, by the triangle inequality, we arrive at
\begin{align}\label{tri}\W_{\vv}(\mu_t^1,\mu_t^2)&\leq \W_\vv(\L_{X_{0,t}^{\gamma^1,\mu^1}}, \L_{X_{0,t}^{\gamma^2,\mu^1}})+ \W_\vv(\L_{X_{0,t}^{\gamma^2,\mu^1}},\L_{X_{0,t}^{\gamma^2,\mu^2}}), \ \ \vv\in(0,1].
\end{align}
and
\begin{align}\label{tri13}\|\mu_t^1-\mu_t^2\|_{var}&\leq \|\L_{X_{0,t}^{\gamma^1,\mu^1}}- \L_{X_{0,t}^{\gamma^2,\mu^1}}\|_{var}+ \|\L_{X_{0,t}^{\gamma^2,\mu^1}}-\L_{X_{0,t}^{\gamma^2,\mu^2}}\|_{var}.
\end{align}
(i)
Let $P_t^{\mu^i}f(x)=\E f(X_t^{x,\mu^i}), x\in\R^d,i=1,2$. Then it follows from \cite[(2.4)]{FYW3} that
\begin{align*}|\nabla P_t^{\mu^i}f|\leq \frac{c}{\sqrt{t}}\|f\|_\infty,\ \ t\in (0,T], f\in\scr B_b(\R^d).
\end{align*}
which implies
\begin{equation*}\begin{split}|P_t^{\mu^i}f(x)-P_t^{\mu^i}f(y)|
&\leq \frac{c}{\sqrt{t}}|x-y|\|f\|_\infty, \ \ t\in (0,T], f\in\scr B_b(\R^d).
\end{split}\end{equation*}
Then there exists a constant $c_1>0$ such that
  \begin{align}\label{kar}\|\L_{X_{0,t}^{\gamma^1,\mu^1}}-\L_{X_{0,t}^{\gamma^2,\mu^1}}\|_{var}\leq  ct^{-\frac{1}{2}}\W_1(\gamma^1,\gamma^2),\ \ t\in (0,T].
\end{align}
This together with \eqref{kpz} implies for any $\vv\in(0,1]$,
  \begin{align}\label{var}\W_{\vv}(\L_{X_{0,t}^{\gamma^1,\mu^1}},\L_{X_{0,t}^{\gamma^2,\mu^1}})\leq  c_2t^{\frac{-1+\vv}{2}}\W_1(\gamma^1,\gamma^2),\ \ t\in (0,T].
\end{align}
Combining \eqref{kar} with \eqref{var}, we conclude that
\begin{align}\label{imp}&\nonumber\sup_{t\in[0,T]}\e^{-\lambda t}\left(t^{\frac{1-\vv}{2}}\W_{\vv}(\L_{X_{0,t}^{\gamma^1,\mu^1}},\L_{X_{0,t}^{\gamma^2,\mu^1}}) +t^{\frac{1}{2}}\|\L_{X_{0,t}^{\gamma^1,\mu^1}}-\L_{X_{0,t}^{\gamma^2,\mu^1}}\|_{var}\right)\\
&\leq (c+c_2)\W_1(\gamma^1,\gamma^2).
\end{align}
(ii) By the boundedness of $b$ and $\sigma$ due to (A1)-(A2), it is easy to see
 $$(P_t^\ast\gamma)(|\cdot|)\leq c(1+\gamma(|\cdot|)),\ \ t\in[0,T],\gamma\in\scr P_1$$
 for some constant $c>0$.
 We derive from \eqref{cty13} for $\gamma=\gamma^2,\nu^i=\mu^i$ that
\begin{align*}\nonumber&\sup_{t\in[0,T]}\e^{-\lambda t}\left(t^{\frac{1-\vv}{2}}\W_{\vv}(\L_{{X}_{0,t}^{\gamma^2,\mu^1}}, \L_{{X}_{0,t}^{\gamma^2,\mu^2}})+t^{\frac{1}{2}}\|\L_{{X}_{0,t}^{\gamma^2,\mu^1}}- \L_{{X}_{0,t}^{\gamma^2,\mu^2}}\|_{var}\right)\\
&\leq c(\eta)\lambda ^{-\frac{\beta\wedge \eta}{2}}\sup_{t\in[0,T]}\e^{-\lambda t}(t^{\frac{1-\eta}{2}}\W_{\eta}(\mu^1_t,\mu^2_t) +t^{\frac{1}{2}}\|\mu^1_t-\mu^2_t\|_{var}).
\end{align*}
Combining this with \eqref{imp} and \eqref{tri}-\eqref{tri13}, for any $\lambda\geq \lambda_0$, we obtain
\begin{equation}\begin{split}\label{LLN}
&\sup_{t\in[0,T]}\e^{-\lambda t}\left(t^{\frac{1-\vv}{2}}\W_{\vv}(\mu_t^1,\mu_t^2)+t^{\frac{1}{2}}\|\mu_t^1-\mu_t^2\|_{var}\right)\\
&\leq \sup_{t\in[0,T]}\e^{-\lambda t}\left(t^{\frac{1-\vv}{2}}\W_{\vv}(\L_{{X}_{0,t}^{\gamma^2,\mu^1}}, \L_{{X}_{0,t}^{\gamma^2,\mu^2}})+t^{\frac{1}{2}}\|\L_{{X}_{0,t}^{\gamma^2,\mu^1}}- \L_{{X}_{0,t}^{\gamma^2,\mu^2}}\|_{var}\right)\\
&+\sup_{t\in[0,T]}\e^{-\lambda t}\left(t^{\frac{1-\vv}{2}}\W_{\vv}(\L_{X_{0,t}^{\gamma^1,\mu^1}},\L_{X_{0,t}^{\gamma^2,\mu^1}})+t^{\frac{1}{2}}\|\L_{X_{0,t}^{\gamma^1,\mu^1}}-\L_{X_{0,t}^{\gamma^2,\mu^1}}\|_{var}\right)\\
&\leq c(\eta)\lambda ^{-\frac{\beta\wedge \eta}{2}}\sup_{t\in[0,T]}\e^{-\lambda t}(t^{\frac{1-\eta}{2}}\W_{\eta}(\mu^1_t,\mu^2_t) +t^{\frac{1}{2}}\|\mu^1_t-\mu^2_t\|_{var})\\
&+(c+c_2)\W_1(\gamma^1,\gamma^2).
\end{split}\end{equation}
Recalling $\lambda_1=\lambda_0\vee(2c(\eta))^{\frac{2}{\beta\wedge\eta}}$, we derive from \eqref{LLN} with $\vv=\eta$ that
\begin{align*}&\sup_{t\in[0,T]}\e^{-\lambda_1 t}\left(t^{\frac{1-\eta}{2}}\W_{\eta}(\mu_t^1,\mu_t^2)+t^{\frac{1}{2}}\|\mu_t^1-\mu_t^2\|_{var}\right)\leq2(c+c_2)\W_1(\gamma^1,\gamma^2).
\end{align*}
This together with \eqref{LLN} implies for any $\vv\in(0,1)$, \begin{align*}&\sup_{t\in[0,T]}\e^{-\lambda_1 t}\left(t^{\frac{1-\vv}{2}}\W_{\vv}(\mu_t^1,\mu_t^2)+t^{\frac{1}{2}}\|\mu_t^1-\mu_t^2\|_{var}\right)\leq 2(c+c_2)\W_1(\gamma^1,\gamma^2).
\end{align*}
The proof is completed.
\end{proof}

\subsection{Proof of Theorem \ref{EEU}(ii)}
The proof of well-posedness is more or less the same as that in \cite[Theorem 1.1(1)]{HWDCDS} since the global integrability condition $|\nabla \sigma|\in L_p^q(T)$ with some $\frac{d}{p}+\frac{2}{q}<1$ in \cite[Theorem 1.1(1)]{HWDCDS}can be replaced by the localized integrability condition $|\nabla \sigma|\in \tilde{L}_p^q(T)$.  Moreover, by the boundedness of $\sigma$, $b^{(0)}$ and (A4), it is standard to derive \eqref{tyq} by It\^{o}'s formula. So, we only verify \eqref{REG}.

Recall $\mu_t^j=P_t^\ast\gamma^j, \ \ j=1,2, t\in [0,T]$ and $P_t^{\mu^j}f(x)=\E f(X_t^{x,\mu^j}), x\in\R^d,j=1,2$.
Let $C_b^1(\R^d)$ be the class of bounded and continuous functions on $\R^d$ with bounded and continuous derivative.
 By an approximation technique, we conclude
\begin{align}\label{dista}\nonumber&\|\mu-\nu\|_{var}=\sup_{f\in C_b^1(\R^d), \|f\|_\infty\leq 1}|\mu(f)-\nu(f)|,\\
&\W_1(\mu,\nu)=\sup_{f\in C_b^1(\R^d), \|\nabla f\|_\infty\leq 1}|\mu(f)-\nu(f)|.
\end{align}
Then it follows from (A1), (A2), (A4) and Theorem \ref{gracd} that we can find a constant $c>0$ depending on $\|b^{(0)}\|_\infty$, $\sigma\sigma^\ast$, $b^{(1)}$ and increasing in $\|b^{(0)}\|_\infty$ such that
\begin{align}\label{gram1}\|\nabla P_t^{\mu^1} f\|_\infty\leq ct^{-\frac{1-i}{2}}\|\nabla ^if\|_\infty, \ \ t\in(0,T], i=0,1,
\end{align}
and
\begin{align}\label{comp1}| P_t^{\mu^1}f(x)- P_t^{\mu^2}f(x)|\leq c\int_{0}^t\|\mu^1_s-\mu^2_s\|_{var}(t-s)^{-\frac{1-i}{2}}\d s\| \nabla^if\|_\infty,\ \ t\in(0,T], i=0,1.
\end{align}
As a result, \eqref{dista} and \eqref{gram1} yield
\begin{align*}\|\L_{{X}_{0,t}^{\gamma^1,\mu^1}}- \L_{{X}_{0,t}^{\gamma^2,\mu^1}}\|_{var}\leq ct^{-\frac{1}{2}}\W_1(\gamma^1,\gamma^2),\ \ \W_1(\L_{{X}_{0,t}^{\gamma^1,\mu^1}}, \L_{{X}_{0,t}^{\gamma^2,\mu^1}})\leq c\W_1(\gamma^1,\gamma^2),
\end{align*}
while \eqref{dista} and \eqref{comp1} imply
\begin{align*}\nonumber& \|\L_{{X}_{0,t}^{\gamma^2,\mu^1}}- \L_{{X}_{0,t}^{\gamma^2,\mu^2}}\|_{var}\leq c\int_{0}^t\|\mu_s^1-\mu_s^2\|_{var}(t-s)^{-\frac{1}{2}}\d s,
\end{align*}
and
\begin{align*}\nonumber& \W_1(\L_{{X}_{0,t}^{\gamma^2,\mu^1}},\L_{{X}_{0,t}^{\gamma^2,\mu^2}})\leq c\int_{0}^t\|\mu_s^1-\mu_s^2\|_{var}\d s.
\end{align*}
By the triangle inequality, we obtain
\begin{align}\label{12v}
\|\mu_t^1-\mu_t^2\|_{var}\leq ct^{-\frac{1}{2}}\W_1(\gamma^1,\gamma^2)+c\int_{0}^t\|\mu_s^1-\mu_s^2\|_{var}(t-s)^{-\frac{1}{2}}\d s,
\end{align}
and
\begin{align}\label{12W}
\W_1(\mu_t^1,\mu_t^2)\leq c\W_1(\gamma^1,\gamma^2)+c\int_{0}^t\|\mu_s^1-\mu_s^2\|_{var}\d s.
\end{align}
Combining \eqref{12v} with \eqref{FKG13} for $\beta=1$, we obtain
\begin{equation}\begin{split}\label{LLN1}
&\sup_{t\in[0,T]}\e^{-\lambda t}t^{\frac{1}{2}}\|\mu_t^1-\mu_t^2\|_{var}\\
&\leq c\W_1(\gamma^1,\gamma^2)+c\sup_{t\in[0,T]}\e^{-\lambda t}t^{\frac{1}{2}}\|\mu_t^1-\mu_t^2\|_{var}\sup_{t\in[0,T]}t^{\frac{1}{2}}\int_{0}^ts^{-\frac{1}{2}}\e^{-\lambda (t-s)}(t-s)^{-\frac{1}{2}}\d s\\
&\leq c\W_1(\gamma^1,\gamma^2)+2c\Gamma(\frac{1}{2})\lambda^{-\frac{1}{2}}\sup_{t\in[0,T]}\e^{-\lambda t}t^{\frac{1}{2}}\|\mu_t^1-\mu_t^2\|_{var}.
\end{split}\end{equation}
Taking $\lambda =16c^2\Gamma(\frac{1}{2})^2$ we derive from \eqref{LLN1} that
\begin{align}\label{l12}&\sup_{t\in[0,T]}\e^{-\lambda_1 t}t^{\frac{1}{2}}\|\mu_t^1-\mu_t^2\|_{var}\leq2c\W_1(\gamma^1,\gamma^2).
\end{align}
Substituting this into \eqref{12W}, we conclude
$$\W_1(\mu_t^1,\mu_t^2)\leq \tilde{c}\W_1(\gamma^1,\gamma^2).$$
Finally, this together with \eqref{l12} asn \eqref{kpz} implies for any $\vv\in(0,1)$, \begin{align*}&\W_{\vv}(\mu_t^1,\mu_t^2)\leq \bar{c}t^{-\frac{1-\vv}{2}}\W_1(\gamma^1,\gamma^2).
\end{align*}
The proof is completed.

\section{Proof of Theorem \ref{T0}}
For $\mu\in \scr P$, consider the time-homogeneous SDEs with parameter $\mu$:
\beq\label{DC} \d {\hat{X}}_t^\mu=b^{(1)}({\hat{X}}_t^\mu)\d t+b^{(0)}({\hat{X}}_t^\mu,\mu)\d t+\si({\hat{X}}_t^\mu)\d W_t,\ \ \ t\geq 0.\end{equation}
According to  \cite[Theorem 2.1]{WSB}, under {\bf(C)}, \eqref{DC} enjoys a unique invariant probability measure denoted by
$\Phi(\mu)$. Let
 $\hat{P}_t^{\mu*}\nu=\L_{\hat{X}_t^\mu}$ for $\hat{X}_t^\mu$ solving \eqref{DC} with initial distribution $\nu\in \scr P.$
Let  $\hat{P}_t^\mu$ be the diffusion semigroup associated   to \eqref{DC}.

Noting $\|b^{(0)}\|_\infty\leq \eta$ and {\bf (C)}, repeating the proof of \cite[Lemma 2.3(2)-(3)]{HW25}, we have the following assertion.
 \beg{lem}\label{L1} Assume {\bf (C)}. Then   the following assertions hold.
\begin{enumerate}
 \item[$(1)$] For any $\theta\in (0,K_2)$,
 $$\Phi(\mu)\big(\e^{\theta|\cdot|^2}\big)<\infty,\ \  \ \mu\in \scr P.$$
\item[$(2)$] There exists a constant
$ \ll(\eta)\in (0,\infty)$ decreasing in $ \eta>0$ such that
\begin{align*}
\W_1({\hat{P}}_t^{\mu *}  \gamma,\hat{P}_t^{\mu*}\tt\gamma)\leq \ff 1 {\ll(\eta)}  \e^{-\lambda(\eta) t}\W_1(\gamma,\tt\gamma),\ \ t\ge 0,\  \gamma,\tt\gamma\in \scr P_{1}.
\end{align*}
Consequently,
\begin{align}\label{gram35}
\|\nn \hat{P}_t^\mu f\|_\infty\le \ff 1 {\ll(\eta)}  \e^{-\lambda(\eta) t}\|\nn f\|_\infty,\ \ \  t\ge 0,\ \mu\in \scr P, f\in C_b^1(\R^d).
\end{align}
\end{enumerate}
 \end{lem}
Next, by {\bf (C)}, \eqref{gram35}, Theorem \ref{gracd} and the semigroup property, we can find a constant $c_0(\eta)\geq 1$ depending on $\eta$, $\sigma\sigma^\ast$, $b^{(1)}$ and increasing in $\eta$ such that
\begin{align}\label{gram13}\|\nabla \hat{P}_t^{\mu} f\|_\infty\leq c_0(\eta)\frac{1}{\lambda (\eta)}\e^{-\lambda(\eta)(t-1)}(t\wedge1)^{-\frac{1}{2}}\|f\|_\infty, \ \ t\geq0, f\in C_b^1(\R^d),
\end{align}
and also find a constant $c_1(T,\eta)>0$ depending on $T,\eta$, $\sigma\sigma^\ast$, $b^{(1)}$ and increasing in $T$ and $\eta$ such that
\begin{align}\label{comp13}| \hat{P}_T^{\mu}f(x)- \hat{P}_T^{\nu}f(x)|\leq \eta c_1(T,\eta)\|\mu-\nu\|_{var}\int_{0}^Ts^{-\frac{1-i}{2}}\d s\| \nabla^if\|_\infty,\ \ i=0,1, f\in C_b^1(\R^d).
\end{align}
With \eqref{gram35}-\eqref{comp13}, Lemma \ref{L1} and \eqref{REG} for $\vv=0$ in hand, we are intending to prove Theorem \ref{T0}. The idea of the proof is more or less the same with that of \cite[Theorem 2.1]{HW25} by replacing $\scr P_{k*}$ therein with $\scr P$.

 \beg{proof}[Proof of Theorem \ref{T0}]
(1) Existence and uniqueness of $\mu_\infty$. By Lemma \ref{L1}(1), we have
 $\Phi: \scr P_1\to  \scr P_1.$
It suffices to show that $\Phi$ is contractive on  $\scr P_1$ under the complete metric
$$\W(\mu,\nu):= \|\mu-\nu\|_{var}+\W_1(\mu,\nu).$$
% By \eqref{gram13}, there exists a constant $c_0(\eta)\in [1,\infty)$ increasing in $\eta$ such that
%\begin{align}\label{kgy}\|P_1^{\mu*}\mu- P_1^{\mu*}\tt\mu\|_{var}\le c_0(\eta) \W_1(\mu,\tt\mu),\ \ \mu,\tilde{\mu}\in \scr P_1.
%\end{align}
Since $\ll(\eta)$ is decreasing in $\eta$, there exists a unique  $T_\eta\in [1,\infty)$ which is increasing in $\eta$ such that
\begin{equation*}  \ff {1+c_0(\eta)} {\ll(\eta)}  \e^{-\lambda(\eta) (T_\eta-1)} =\ff 1 2,
\end{equation*}   so that \eqref{gram35}-\eqref{gram13} and the fact  $\int_{\R^d}f(x)(\hat{P}_t^{\mu*}\gamma)(\d x) =\int_{\R^d} \hat{P}_t^{\mu}f(x)\gamma(\d x), f\in\scr B_b(\R^d)$ imply
\begin{equation}\label{s-1}\W(\hat{P}_{T_\eta}^{\mu*}\gamma, \hat{P}_{T_\eta}^{\mu*}\tt\gamma )=\|\hat{P}_{T_\eta}^{\mu*}\gamma- \hat{P}_{T_\eta}^{\mu*}\tt\gamma\|_{var} + \W_1(\hat{P}_{T_\eta}^{\mu*}\gamma, \hat{P}_{T_\eta}^{\mu*}\tt\gamma ) \le \ff 1 2 \W_1(\gamma,\tt\gamma),\ \ \gamma,\tt\gamma\in \scr P_1.
\end{equation}
Next, by \eqref{comp13} and \eqref{dista}, we obtain
 \beq\label{14} \beg{split}  & \|\hat{P}_{T_\eta}^{\mu*} \Phi(\nu)- \hat{P}_{T_\eta}^{\nu*} \Phi(\nu)\|_{var}=\sup_{f\in C_b^1(\R^d),\|f\|_{\infty}\le 1} \big|\Phi(\nu)(\hat{P}_{T_\eta}^\mu f- \hat{P}_{T_\eta}^\nu f)\big| \\
 &\le\eta \|\mu-\nu\|_{var} \int_0^{T_\eta} c_1(T_\eta,\eta) t^{-\ff 1 2}\d t= 2\eta c_1(T_\eta,\eta)T_\eta^{1/2}\|\mu-\nu\|_{var},\end{split}\end{equation}
and
\beg{align*} & \W_1(\hat{P}_{T_\eta}^{\mu*} \Phi(\nu), \hat{P}_{T_\eta}^{\nu*} \Phi(\nu)) =\sup_{f\in C_b^1(\R^d), \|\nn f\|_\infty\le 1} \big|\Phi(\nu)(\hat{P}_{T_\eta}^\mu f- \hat{P}_{T_\eta}^\nu f)\big|  \le  \eta c_1(T_\eta,\eta)T_\eta  \|\mu-\nu\|_{var}.
\end{align*}
Combining this with \eqref{14}, we derive
$$\W(\hat{P}_{T_\eta}^{\mu*} \Phi(\nu), \hat{P}_{T_\eta}^{\nu*} \Phi(\nu))\le \eta c_1(T_\eta,\eta)\max\Big\{2T_\eta^{1/2}, T_\eta\Big\}\W(\mu,\nu).$$
This together with \eqref{s-1} for $\gamma=\Phi(\mu)$ and $\tilde{\gamma}=\Phi(\nu)$ and $\hat{P}_t^{\mu*} \Phi(\mu)=\Phi(\mu), \ \ t\ge 0,\ \ \mu\in \scr P$ implies
\beq\label{13} \beg{split} &\W(\Phi(\mu),\Phi(\nu))= \W(\hat{P}_{T_\eta}^{\mu*}\Phi(\mu), \hat{P}_{T_\eta}^{\nu*}\Phi(\nu))\\
&\le  \W(\hat{P}_{T_\eta}^{\mu*}\Phi(\mu), \hat{P}_{T_\eta}^{\mu*}\Phi(\nu))+  \W(\hat{P}_{T_\eta}^{\mu*}\Phi(\nu), \hat{P}_{T_\eta}^{\nu*}\Phi(\nu))\\
&\le \ff 1 2\W(\Phi(\mu),\Phi(\nu))+\eta c_1(T_\eta,\eta)\max\Big\{2T_\eta^{1/2}, T_\eta\Big\}\W(\mu,\nu),\ \ \ \mu,\nu\in \scr P_1.\end{split}\end{equation}
Note that there exists a constant $\eta_0>0$ depending only on $(a,b^{(1)})$ such that
 $$\eta c_1(T_\eta,\eta)\max\Big\{2T_\eta^{1/2}, T_\eta\Big\}\le \ff 1 4,\ \ \eta\in (0,\eta_0].$$
Substituting this into \eqref{13}, we conclude that $\Phi$ is $\W$-contractive when $\eta\le \eta_0.$

 (2)  $\W_1$-exponential ergodicity. For $\mu\in \scr P$,
 consider  the decoupled time-inhomogeneous SDE
\begin{equation}\label{inh} \d \bar X_{s,t}^\mu=b^{(1)}(\bar X_{s,t}^\mu)\d t+b^{(0)}(\bar X_{s,t}^\mu, P_t^*\mu)\d t +\si(\bar X_{s,t}^\mu) \d W_t,\ \ t\ge s\geq 0.
\end{equation}
Under {\bf (C)}, this SDE is well-posed. For any $\nu\in \scr P$, let
$\bar P_{s,t}^{\mu *}\nu:=\L_{\bar X_{s,t}^\mu}$ for $\bar X_{s,t}^\mu$ solving \eqref{inh} with $\L_{\bar X_{s,s}^\mu}=\nu.$
Let $\bar P_{s,t}^\mu $ be the associated semigroup to \eqref{inh}.
 Simply denote $\bar P_{t}^\mu=\bar P_{0,t}^\mu$ and $\bar P_{t}^{\mu *}=\bar P_{0,t}^{\mu *}$.
Then
\beq\label{LX} P_t^*\mu= \bar P_t^{\mu*}\mu,\ \ \ \mu_\infty= P_t^*\mu_\infty= \bar P_t^{\mu_\infty*}\mu_\infty= \hat{P}_t^{\mu_\infty*}\mu_\infty.\end{equation}
 By \eqref{LX} and Lemma \ref{L1}(2), it holds
 \beq\label{EN1} \beg{split}&\W_1(P_t^*\mu,\mu_\infty)=  \W_1(\bar P_t^{\mu*}\mu,\hat{P}_t^{\mu_\infty*}\mu_\infty)\\
 &\le  \W_1(\bar P_t^{\mu*}\mu, \hat{P}_t^{\mu_\infty*}\mu)+  \W_1(\hat{P}_t^{\mu_\infty*}\mu, \hat{P}_t^{\mu_\infty*}\mu_\infty)\\
 &\le \W_1(\bar P_t^{\mu*}\mu, \hat{P}_t^{\mu_\infty*}\mu)+ \ff 1 {\ll(\eta)}  \e^{-\ll(\eta)t}\W_1(\mu,\mu_\infty),\ \ t\ge 0,\ \ \mu\in\scr P_1.\end{split}\end{equation}
Let  $\tilde{T}_\eta\in (0,\infty)$ be  increasing in $\eta$  such that
 $\ff 1 {\ll(\eta)} \e^{-\ll(\eta)\tilde{T}_\eta}=\ff 1 4.$ Then \eqref{EN1} implies
 \beq\label{EN135} \beg{split}&\W_1(P_t^*\mu,\mu_\infty)\le \W_1(\bar P_t^{\mu*}\mu, \hat{P}_t^{\mu_\infty*}\mu)+ \frac{1}{4}\W_1(\mu,\mu_\infty) \ \ t\ge \tilde{T}_\eta,\ \ \mu\in\scr P_1.\end{split}\end{equation}
Again by  {\bf (C)}, Theorem \ref{gracd} for $i=1$, we find a constant $c_2(t,\eta)>0$ depending on $t,\eta$, $\sigma\sigma^\ast$, $b^{(1)}$ and increasing in $t$ and $\eta$ such that
\begin{align}\label{com13}
\nonumber \W_1(\bar P_t^{\mu*} \mu, \hat{P}_t^{\mu_\infty*} \mu)&= \sup_{f\in C_b^1(\R^d),\|\nn f\|_\infty\le 1} |\mu(\bar P_t^\mu f-\hat{P}_t^{\mu_\infty}f)|\\
 &\leq  \eta c_2(t,\eta) \int_0^t \|P_s^*\mu-\mu_\infty\|_{var} \d s,\ \ t\geq 0.\end{align}
 On the other hand, by \eqref{REG} for $\vv=0$ and $T=(2\tilde{T}_\eta)\vee1$, we have
\beq\label{PLN} \|P_{t}^*\mu-\mu_\infty\|_{var}= \|P_{t}^*\mu-P_t^{*} \mu_\infty\|_{var}\le \kk_1(\eta) t^{-\ff 1 2}\W_1(\mu,\mu_\infty),\ \ \  t\in (0, (2\tilde{T}_\eta)\vee 1]\end{equation}
for some  constant $\kk_1(\eta)\in (0,\infty)$  increasing in $\eta$.
Combining \eqref{PLN} with \eqref{com13}, we find a constant $\kk_2(\eta)\in (0,\infty)$ increasing in $\eta$ such that
 \beq\label{N3} \W_1(\bar P_t^{\mu*} \mu, \hat{P}_t^{\mu_\infty*} \mu) \le  \eta \kk_2(\eta) \W_1(\mu,\mu_\infty),\ \ t\in [0,2 \tilde{T}_\eta].\end{equation}
 Combining \eqref{N3} with \eqref{EN135} and taking   $\eta_1\leq \eta_0$ such that  {$\eta_1 \kk_2(\eta_1)\le \ff 1 4$},  we deduce that $\eta\le \eta_1$ implies
 $$\W_1(P_t^*\mu,\mu_\infty)\le \ff 1 2 {\W_1(\mu,\mu_\infty),\ \ \ t\in [\tilde{T}_\eta, 2\tilde{T}_{\eta}]}.$$
This together with \eqref{EN1}, \eqref{N3} and the semigroup property of $P_t^\ast$ implies \eqref{EX1} for some constants $c,\ll>0$.
Combining \eqref{EX1}  with \eqref{PLN} and the semigroup property of $P_t^*$, we derive
 \beg{align*}&\|P_t^*\gg-\mu_\infty\|_{var}=\|P_{1}^*P_{t-1}^*\gg-\mu_\infty\|_{var}\le \kk_1(\eta)   c\e^{-\ll (t-1)} \W_1(\gg,\mu_\infty),\ \ \ t>1.\end{align*}
Combining this with \eqref{PLN} for $t\le 1$, \eqref{EX2} holds for the same $\ll>0$ but a  different constant  $c >0$.
   \end{proof}
\section{SDEs with monotone drift plus a bounded measurable perturbation}
 Consider two different classical SDEs:
\begin{align}\label{X10}\d X_{s,t}=b_t^{(1)}(X_{s,t})\d t+b^{(0)}_t(X_{s,t})\d t+\sigma_t(X_{s,t})\d W_t,\ \ t\geq s\geq 0,
\end{align}
and
\begin{align}\label{X1t}
\d \tilde{X}_{s,t}=b_t^{(1)}(\tilde{X}_{s,t})\d t+\tilde{b}^{(0)}_t(\tilde{X}_{s,t})\d t+\sigma_t(\tilde{X}_{s,t})\d W_t, \ \ t\geq s\geq 0.
\end{align}
where $b^{(1)}, b^{(0)}, \tilde{b}^{(0)}:[0,\infty)\times \R^d\to\R^d$, $\si: [0,\infty)\times\R^d\to \R^d\otimes\R^m$ are measurable and locally bounded. In this part, we will investigate the gradient estimate for \eqref{X10} and \eqref{X1t} and derive a comparison estimate between them. To this end, $b^{(1)}$ is assumed to satisfy monotone condition while $b^{(0)}, \tilde{b}^{(0)}$ are only bounded measurable.
More precisely, we make the following assumption.
\begin{enumerate}
\item[{\bf (D)}] $b^{(0)}$ and $\tilde{b}^{(0)}$ are bounded measurable. There exist constants $\delta_1\geq \delta_2>0$ such that
  $$ \delta_2\leq\sigma\sigma^\ast\leq \delta_1.$$
 Moreover, there exists a constant  $K_1\geq 0$ such that
      $$\<b_t^{(1)}(x)-b_t^{(1)}(y),x-y\>\leq K_1|x-y|^2,$$
      and $$\|\sigma_t(x)-\sigma_t(y)\|_{HS}^2\leq K_1|x-y|^2.$$
      \end{enumerate}
Under {\bf(D)}, \eqref{X10} and \eqref{X1t} are well-posed and non-explosive. Let $P_{s,t}$ and $\tilde{P}_{s,t}$ be the associated semigroups to \eqref{X10} and \eqref{X1t} respectively, i.e.
$$P_{s,t} f(x)=\E f(X_{s,t}^x),\ \ \tilde{P}_{s,t} f(x)=\E f(\tilde{X}_{s,t}^x),\ \  f\in\scr B_b(\R^d),t\geq s\geq 0, x\in\R^d,$$
where $X_{s,t}^x$($\tilde{X}_{s,t}^x$) solve \eqref{X10}(\eqref{X1t}) with $X_{s,s}^x=x\in\R^d$($\tilde{X}_{s,s}^x=x\in\R^d$). For simplicity, we denote $P_s=P_{0,s}$($\tilde{P}_s=\tilde{P}_{0,s}$).
The main result in this part is stated as follows.
\begin{thm}\label{gracd} Assume {\bf (D)}. Then for any $T>0$, there exists a constant $c>0$ depending on $T,\|b^{(0)}\|_\infty, \|\tilde{b}^{(0)}\|_\infty $, $\delta_1,\delta_2, K_1$ and increasing in $T$ and $\|b^{(0)}\|_\infty+ \|\tilde{b}^{(0)}\|_\infty $ such that for any $f\in C_b^1(\R^d)$,
\begin{align}\label{gratn}\|\nabla P_{s,t} f\|_\infty+\|\nabla \tilde{P}_{s,t}f\|_\infty\leq c(t-s)^{-\frac{1-i}{2}}\|\nabla ^if\|_\infty, \ \ 0\leq s<t\leq T, i=0,1,
\end{align}
and
\begin{align}\label{compa}| \tilde{P}_tf(x)- P_tf(x)|\leq c\int_{0}^t\|\tilde{b}^{(0)}_s-b_s^{(0)}\|_{\infty}(t-s)^{-\frac{1-i}{2}}\d s\| \nabla^if\|_\infty,\ \ t\in(0,T], i=0,1.
\end{align}
\end{thm}
\begin{rem} When $b^{(0)}$, $b^{(1)}$, $\sigma$ are independent of $t$, {\bf (D)} implies
\begin{align*}&\|\sigma(x)-\sigma(y)\|_{HS}^2+2\<b^{(1)}(x)-b^{(1)}(y)+b^{(0)}(x)-b^{(0)}(y),x-y\>\\
&\leq 3K_1|x-y|^2+4\|b^{(0)}\|_\infty|x-y|.
\end{align*}
If in addition $b^{(0)}$, $b^{(1)}$ are continuous, then according to \cite[Theorem 3.4(a)]{PR}, it holds
$$\|\nabla P_t f\|\leq c_t\|f\|_\infty,\ \ f\in \scr B_b(\R^d)$$
for $$c_t=\inf_{r>0}\left\{\frac{\int_0^r\exp\{\frac{1}{4\lambda_0}\int_0^sg(u)\d u\}\d s}{2\lambda_0 t}+\frac{1}{\int_0^r\exp\{-\frac{1}{4\lambda_0}\int_0^sg(u)\d u\}\d s}\right\}$$
and $g(r)=3K_1 r+4\|b^{(0)}\|_\infty$. We should remark that the proof of Theorem \ref{gracd} is quite different from that of  \cite[Theorem 3.4(a)]{PR}, where a coupling method is adopted.

\end{rem}

\begin{rem} If in addition, $b^{(1)}_t$ is Lipschitz continuous, i.e. there exists a cosntant $K\geq 0$ such that $$|b_t^{(1)}(x)-b_t^{(1)}y)|\leq K|x-y|,$$
then the assertion in Theorem \ref{gracd} is known, see for instance \cite{HRW25}. However, in this case, it turns out that the constant $c$ in \eqref{gratn} depends on the Lipschitz constant $K$.
\end{rem}
We first prove Theorem \ref{gracd} under the assumption that $b_t^{(0)}$ and $\tilde{b}_t^{(0)}$ are Lipschitz continuous.
%(1) Without loss of generality, we can start from deterministic initial value. We first assume that $b^{(1)}$ is Lipschitz and uniformly dissipative. Using the regularization approximation on $b^{(0)}$, we can remove the uniform continuous condition on $b^{(0)}$. Here, we may use the Krylov estimate.

%(2) By Duhamel's formula, see \cite[(5.33)]{HRW25}, we have
% \begin{equation*} \bar P_{s,t}^\mu =\hat P_{t-s} f+ \int_s^t \bar P_{s,r}^\mu \<b^{(0)}(\cdot, P_r^*\mu), \nn\hat P_{t-s}f\> \d r,\ \ t\ge s\ge 0,\ f\in\B_b(\R^d).
% \end{equation*}

%Using Yosida's approximation, and the Duhamel formula
%$$P_tf=P_t^n f+\int_0^tP_{s}\<b^{(1)}-b^{(1),n}, \nabla P^n_{t-s}f\>$$
%Using the gradient in Priola and Wang, the dominated convergence theorem.

\begin{lem}\label{gracd1} Assume {\bf (D)} and $b_t^{(0)}$ and $\tilde{b}_t^{(0)}$ are Lipschitz continuous with Lipschitz constant independent of $t$. Then the assertions in Theorem \ref{gracd} hold.
\end{lem}
\begin{proof} It is sufficient to prove \eqref{gratn} for $s=0$. We divide the proof into two steps.

{\bf Step (i)} Assume $b_t^{(1)}$, $b_t^{(0)}$ and $\tilde{b}_t^{(0)}$ are Lipschitz continuous with Lipschitz constant independent of $t$. Denote $[f]_{Lip}$ the Lipschitz constant of $f$ and $[b^{(1)}]_{Lip,\infty}=\sup_{t\geq 0}[b_t^{(1)}]_{Lip}$. Consider the reference SDE
\begin{align}\label{reSDE}\d \hat{X}_{s,t}^1=b^{(1)}_{t}(\hat{X}_{s,t}^1)\d t+\sigma_t(\hat{X}_{s,t}^1)\d W_t.
\end{align}
Let $\hat{P}_{s,t}$ be the associated semigroup to \eqref{reSDE} and we simply denote $\hat{P}_{s}=\hat{P}_{0,s}$.
Under {\bf (D)}, we can find a constant ${\bf c_G}>0$ depending on $T,\delta_1,\delta_2, K_1$ (independent of the Lipschitz constant of $b^{(1)}_t$) and increasing in $T$ such that for any $f\in C_b^1(\R^d)$,
\begin{align}\label{regra}\|\nabla \hat{P}_{s,t} f\|_\infty\leq {\bf c_G}(t-s)^{-\frac{1-i}{2}}\|\nabla ^if\|_\infty,\ \ 0\leq s<t\leq T,i=0,1.
\end{align}
This means
\begin{align}\label{kpinf}\|\nabla \hat{P}_{s,t} \|_{\infty\to\infty}:=\sup_{\|f\|_\infty\leq 1}\|\nabla \hat{P}_{s,t}f\|_\infty \leq {\bf c_G}(t-s)^{-\frac{1}{2}},\ \ 0\leq s<t\leq T.
\end{align}
Since $b_t^{(1)}$ is Lipschitz continuous  with Lipschitz constant independent of $t$ and $b^{(0)}$ is bounded,
we may use Duhamel's formula(see for instance \cite[(5.33)]{HRW25}) to obtain
\begin{align}\label{duham}P_tf=\hat{P}_tf+\int_0^tP_{s}\<b^{(0)}_s, \nabla \hat{P}_{s,t}f\>\d s,\ \ f\in \scr B_b(\R^d),
\end{align}
$$\tilde{P}_tf=\hat{P}_t f+\int_0^t\tilde{P}_{s}\<\tilde{b}^{(0)}_s, \nabla \hat{P}_{s,t}f\>\d s,\ \ f\in \scr B_b(\R^d),$$
and
\begin{align}\label{ptk}\tilde{P}_tf=P_t f+\int_0^t\tilde{P}_{s}\<\tilde{b}^{(0)}_s-b^{(0)}_s, \nabla P_{s,t}f\>\d s,\ \ f\in \scr B_b(\R^d).
\end{align}
Moreover, by \cite[Proposition 5.2]{HRW25}, there exists a constant $C_0>0$ depending on $\delta_1,\delta_2, K_1$, $\|b^{(0)}\|_\infty$ and $[b^{(1)}]_{Lip,\infty}$ such that \eqref{gratn} holds,
which implies
\begin{align}\label{p-p}H_\lambda(T):=\sup_{t\in[0,T]}\e^{-\lambda t}\sqrt{t}\|\nabla P_t\|_{\infty\to\infty}<\infty, \ \ \lambda>0,T>0.
\end{align}
Next, we verify that the constant $C_0$ can be independent of $[b^{(1)}]_{Lip,\infty}$ due to the monotonicity of $b_t^{(1)}$. It is sufficient to prove \eqref{gratn} for $P_t$.
For any $f\in \scr B_b(\R^d)$ with $\|f\|_\infty\leq 1$, we derive from \eqref{duham}, \eqref{p-p}, \eqref{regra} for $i=0$ and the dominated convergence theorem that
\begin{align}\label{ctyk1}
\nabla P_tf=\nabla \hat{P}_t f+\int_0^t\nabla P_{s}\<b^{(0)}_s, \nabla \hat{P}_{s,t}f\>\d s,
\end{align}
which together with \eqref{kpinf} implies
\begin{align*}
\|\nabla P_t\|_{\infty\to\infty}&\leq \|\nabla \hat{P}_t\|_{\infty\to \infty}+\int_0^t\|\nabla P_{s}\|_{\infty\to\infty}\|b^{(0)}\|_\infty\|\nabla \hat{P}_{s,t}\|_{\infty\to\infty}\d s\\
&\leq {\bf c_G}t^{-\frac{1}{2}}+\|b^{(0)}\|_\infty\int_0^t\|\nabla P_{s}\|_{\infty\to\infty}{\bf c_G}(t-s)^{-\frac{1}{2}}\d s.
\end{align*}
This gives
\begin{align}\label{Hla}
H_\lambda (T)&\leq {\bf c_G}+H_\lambda(T)\|b^{(0)}\|_\infty\sup_{t\in[0,T]}\sqrt{t}\int_0^t\e^{-\lambda (t-s)}s^{-\frac{1}{2}}{\bf c_G}(t-s)^{-\frac{1}{2}}\d s.
\end{align}
Note that the FKG inequality implies
\begin{align*}&\sqrt{t}\int_0^t\e^{-\lambda (t-s)}s^{-\frac{1}{2}}{\bf c_G}(t-s)^{-\frac{1}{2}}\d s\\
&\leq {\bf c_G}t^{\frac{3}{2}}\frac{1}{t}\int_0^ts^{-\frac{1}{2}}\d s\frac{1}{t}\int_0^t\e^{-\lambda (t-s)}(t-s)^{-\frac{1}{2}}\d s\\
&\leq 2{\bf c_G}\int_0^\infty\e^{-\lambda s}s^{-\frac{1}{2}}\d s=2\Gamma(\frac{1}{2}){\bf c_G}\lambda ^{-\frac{1}{2}}.
\end{align*}
Taking $\lambda_0=16\|b^{(0)}\|_\infty^2\Gamma(\frac{1}{2})^2{\bf c_G}^2$, we have
$\|b^{(0)}\|_\infty2\Gamma(\frac{1}{2}){\bf c_G}\lambda_0 ^{-\frac{1}{2}}= \frac{1}{2}.$
Combining \eqref{Hla} with \eqref{p-p}, we arrive at
\begin{align*}
\sup_{t\in[0,T]}\e^{-\lambda_0 t}\sqrt{t}\|\nabla P_t\|_{\infty\to\infty}\leq 2{\bf c_G}.
\end{align*}
which implies
\begin{align}\label{suHl3}
\|\nabla P_t\|_{\infty\to\infty}\leq 2{\bf c_G}\e^{\lambda_0 t}t^{-\frac{1}{2}}, \ \ t\in(0,T].
\end{align}
So, for $f\in C_b^1(\R^d)$, it follows from \eqref{ctyk1}, \eqref{regra} and \eqref{suHl3} that
\begin{align*}
\|\nabla P_tf\|_{\infty}&\leq \|\nabla \hat{P}_t f\|_{\infty}+\int_0^t\|\nabla P_{s}\|_{\infty\to\infty}\|b^{(0)}\|_{\infty} \|\nabla \hat{P}_{t-s}f\|_{\infty}\d s\\
&\leq {\bf c_G}\|\nabla f\|_\infty+\int_0^t2{\bf c_G}^2\e^{\lambda_0 s}s^{-\frac{1}{2}}\|b^{(0)}\|_{\infty} \|\nabla f\|_{\infty}\d s\\
& \leq \left({\bf c_G}+4{\bf c_G}^2\e^{\lambda_0 t}\sqrt{t}\|b^{(0)}\|_{\infty}\right)\|\nabla f\|_\infty, \ \ t\in(0,T].
\end{align*}
This and \eqref{suHl3} yield \eqref{gratn}, which combined with \eqref{ptk} implies \eqref{compa} immediately.

{\bf Step (ii)} Assume $b_t^{(0)}$ and $\tilde{b}_t^{(0)}$ are Lipschitz continuous with Lipschitz constant independent of $t$. We use the Yosida approximation for $b^{(1)}$.
Let
    $
        \tilde{b}_t^{(1)}(x):=b_t^{(1)}(x)
        -K_1x,\ \ x\in\R^{d}.
    $
    For any $m\geq 1$, let
    $$
        \tilde{b}^{(1),m}_t(x)
        :=m\left[
        \left(\operatorname{id}
        -\frac{1}{m}\tilde{b}_t^{(1)}\right)^{-1}(x)-x
        \right],\quad x\in\R^d,
    $$
    where $\mathrm{id}$ is the identity map on $\R^d$.
     By \cite[Proposition D.11]{GZ}(see also \cite[Section 2]{GRW}), we have
\begin{align}\label{abc}
|\tilde{b}^{(1),m}|\leq|\tilde{b}^{(1)}|,\ \ m\geq 1, \ \ \lim_{m\to\infty}\tilde{b}^{(1),m}=\tilde{b}^{(1)},
\end{align}
and
$$\<\tilde{b}_t^{(1),m}(x)-\tilde{b}_t^{(1),m}(\tilde{x}),x-\tilde{x}\>\leq 0,\ \  |\tilde{b}^{(1),m}_t(x)-\tilde{b}^{(1),m}_t(\tilde{x})|\leq 2m|x-\tilde{x}|,\ \ t\geq 0, x,\tilde{x}\in \R^d, m\geq 1. $$
Let
    \begin{align*}b_t^{(1),m}(x):=\tilde{b}^{(1),m}_t(x)+K_1x, \ \ t\geq 0,x\in\R^d,m\geq 1.
    \end{align*}
Then one has
\begin{align}\label{bt0}|b_t^{(1),m}(0)|=|\tilde{b}_t^{(1),m}(0)|\leq |\tilde{b}_t^{(1)}(0)|=|b_t^{(1)}(0)|,\ \ t\geq 0,m\geq 1,
\end{align}
\begin{align}\label{bhn}|b^{(1),m}_t(x)-b^{(1),m}_t(\tilde{x})|\leq (2m+K_1)|x-\tilde{x}|,\ \ x,\tilde{x}\in\R^d,t\geq 0,m\geq 1,
\end{align} and
    \begin{equation}\label{jg4dv7}
        \<b_t^{(1),m}(x)-b_t^{(1),m}(y),x-y\>\leq K_1|x-y|^2,\ \ m\geq 1, x,y\in\R^d, t\geq 0.
    \end{equation}
Consider
\begin{align}\label{Xm}
\d X_{s,t}^{m,x}=b_t^{(1),m}(X_{s,t}^{m,x})\d t+b^{(0)}_t(X_{s,t}^{m,x})\d t+\sigma_t(X_{s,t}^{m,x})\d W_t,\ \ t\geq s\geq 0, X_{s,s}^{m,x}=x\in\R^d,
\end{align}
and
\begin{align}\label{tXm}
\d \tilde{X}_{s,t}^{m,x}=b_t^{(1),m}(\tilde{X}_{s,t}^{m,x})\d t+\tilde{b}^{(0)}_t(\tilde{X}_{s,t}^{m,x})\d t+\sigma_t(\tilde{X}_{s,t}^{m,x})\d W_t,\ \ t\geq s\geq 0,\tilde{X}_{s,s}^{m,x}=x\in\R^d.
\end{align}
Let $P_{s,t}^m$ and $\tilde{P}_{s,t}^m$ be the associated semigroups to \eqref{Xm} and \eqref{tXm} respectively.
Then by \eqref{bhn}, \eqref{jg4dv7} and {\bf Step (i)}, there exist a constant $c>0$ depending on $T,\|b^{(0)}\|_\infty, \|\tilde{b}^{(0)}\|_\infty $, $\delta_1,\delta_2, K_1$ and increasing in $T$ and $\|b^{(0)}\|_\infty+ \|\tilde{b}^{(0)}\|_\infty$ such that for any $f\in C_b^1(\R^d)$,
\begin{align}\label{myk}\|\nabla \tilde{P}_{s,t}^m f\|_\infty+\|\nabla P_{s,t}^m f\|_\infty\leq c(t-s)^{-\frac{1-i}{2}}\|\nabla ^if\|_\infty, \ \ 0\leq s<t\leq T, i=0,1,
\end{align}
and
\begin{align}\label{cmz}
\nonumber&| \tilde{P}^m_tf(x)- P_t^mf(x)|\\
&\leq c\int_{0}^t\|\tilde{b}^{(0)}_s-b_s^{(0)}\|_{\infty}(t-s)^{-\frac{1-i}{2}}\d s\| \nabla^if\|_\infty,\ \ t\in(0,T], i=0,1, x\in\R^d.
\end{align}
It follows from \eqref{myk} that
\begin{align}\label{dif}\nonumber&| P_{s,t}^{m}f(x)- P_{s,t}^{m}f(y)|+| \tilde{P}_{s,t}^{m}f(x)- \tilde{P}_{s,t}^{m}f(y)|\\
&\leq c(t-s)^{-\frac{1-i}{2}}\|\nabla ^if\|_\infty|x-y|,\ \  0\leq s<t\leq T, i=0,1,x,y\in\R^d.
\end{align}
Letting $m\to\infty$ in \eqref{cmz} and \eqref{dif}, we derive \eqref{gratn} and \eqref{compa} by Lemma \ref{yap} below and the dominated convergence theorem.
\end{proof}
We are now in the position to complete the proof of Theorem \ref{gracd}.
\begin{proof}[Proof of Theorem \ref{gracd}] To use the result in Lemma \ref{gracd1}, we use the approximation by mollifier.
Let $\rho\in C_0^\infty(\R^{d}), \rho\geq 0$ with $\int_{\R^{d}}\rho(x)\d x=1$, $\rho_n(x)=n^{d}\rho(nx)$. Let $b_t^{(0),n}=\rho_n\ast b_t^{(0)}, \tilde{b}_t^{(0),n}=\rho_n\ast \tilde{b}^{(0)}_t$.
Then for any $t\in[0,T]$, $b_t^{(0),n}$($\tilde{b}_t^{(0),n}$) converges to $b_t^{(0)}$ ($\tilde{b}_t^{(0)}$) almost surely as $n\to\infty$. Moreover, $b_t^{(0),n}, \tilde{b}_t^{(0),n}$ are Lipschitz continuous with Lipschitz constant independent of $t$ and for any $n\geq 1$,
\begin{align}\label{bci}\|b_t^{(0),n}\|_{\infty}\leq \|b^{(0)}_t\|_\infty, \ \ \|\tilde{b}^{(0),n}_t\|_{\infty}\leq \|\tilde{b}^{(0)}_t\|_\infty, \ \ \|b_t^{(0),n}-\tilde{b}^{(0),n}_t\|_\infty\leq \|b_t^{(0)}-\tilde{b}^{(0)}_t\|_\infty.
\end{align}
Consider
\begin{align}\label{Xnt}\d X_{s,t}^{\<n\>,x}=b_t^{(1)}(X_{s,t}^{\<n\>,x})\d t+b_t^{(0),n}(X_{s,t}^{\<n\>,x})\d t+\sigma_t(X_{s,t}^{\<n\>,x})\d W_t,\ \ t\geq s\geq 0,X_{s,s}^{\<n\>,x}=x,
\end{align}
and
\begin{align}\label{cty35}\d \tilde{X}_{s,t}^{\<n\>,x}=b_t^{(1)}(\tilde{X}_{s,t}^{\<n\>,x})\d t+\tilde{b}_t^{(0),n}(\tilde{X}_{s,t}^{\<n\>,x})\d t+\sigma_t(\tilde{X}_{s,t}^{\<n\>,x})\d W_t,\ \ t\geq s\geq 0,\tilde{X}_{s,s}^{\<n\>,x}=x.
\end{align}
 Let $P_{s,t}^{\<n\>}$ and $\tilde{P}_{s,t}^{\<n\>}$ be the associated semigroups to \eqref{Xnt} and \eqref{cty35} respectively. By Lemma \ref{gracd1} and \eqref{bci}, we can find a constant $c>0$ depending on $T ,\|b^{(0)}\|_\infty, \|\tilde{b}^{(0)}\|_\infty $, $\delta_1,\delta_2, K_1$ and increasing in $T$ and $\|b^{(0)}\|_\infty+ \|\tilde{b}^{(0)}\|_\infty$ such that for any $f\in C_b^1(\R^d)$,
$$\|\nabla P_{s,t}^{\<n\>} f\|_\infty+\|\nabla \tilde{P}_{s,t}^{\<n\>} f\|_\infty\leq c(t-s)^{-\frac{1-i}{2}}\|\nabla ^if\|_\infty, \ \ 0\leq s<t\leq T, i=0,1,$$
and
$$| \tilde{P}^{\<n\>}_tf(x)- P_t^{\<n\>}f(x)|\leq c\int_{0}^t\|\tilde{b}^{(0)}_s-b_s^{(0)}\|_{\infty}(t-s)^{-\frac{1-i}{2}}\d s\| \nabla^if\|_\infty,\ \ t\in(0,T], i=0,1, x\in\R^d$$
By the same argument to derive \eqref{gratn} and \eqref{compa} from \eqref{cmz} and \eqref{dif}, it remains to
prove
\begin{align}\label{kct}\lim_{n\to\infty}[|P_{s,t}^{\<n\>}f-P_{s,t} f|+|\tilde{P}_{s,t}^{\<n\>}f-\tilde{P}_{s,t} f|]=0,\ \ f\in C_b^1(\R^d).
\end{align}
It is sufficient to verify \eqref{kct} for $P_{s,t}$ and $P_{s,t}^{\<n\>}$ with $s=0$. Again for simplicity, we denote $X_t^{\<n\>,x}=X_{0,t}^{\<n\>,x}$ and $P_{t}^{\<n\>}=P_{0,t}^{\<n\>}$. Let $s=0$ and
rewrite \eqref{Xnt} as
\begin{align*}\d X_t^{\<n\>,x}=b_t^{(1)}(X_t^{\<n\>,x})\d t+b^{(0)}_t(X_t^{\<n\>,x})\d t+\sigma_t(X_t^{\<n\>,x})\d W^n_t,
\end{align*}
where $$\d W^n_t=\d W_t-\xi_t^n\d t,\ \ \xi_t^n=[\sigma_t(\sigma_t\sigma^\ast_t)^{-1}](X_t^{\<n\>,x})(b^{(0)}_t(X_t^{\<n\>,x})-b^{(0),n}_t(X_t^{\<n\>,x})).$$
By Girsanov's theorem, $(W^n_t)_{t\in[0,T]}$ is a $d$-dimensional Brownian motion under the probability measure $\d \Q_T^n=R_T^n\d \P$ for the martingale
$$R_t^n=\exp\left\{\int_0^t\<\xi_s^n,\d W_s\>-\frac{1}{2}\int_0^t|\xi_s^n|^2\d s\right\},\ \ t\in[0,T].$$
This means
$$P_t f(x)=\E^{\Q_T^n} f(X_t^{\<n\>,x})=\E(R_T^n f(X_t^{\<n\>,x}))=\E(R_t^n f(X_t^{\<n\>,x})), \ \ f\in \scr B_b(\R^d).$$
By Young's inequality, we conclude that
\begin{align*}P_t\log f(x)&\leq \log P_t^{\<n\>} f(x)+\E(R_t^n\log R_t^n)\\
&=\log P_t^{\<n\>} f(x)+\frac{1}{2}\E^{\Q_T^n}\int_0^t|\xi_s^n|^2\d s\\
&=\log P_t^{\<n\>} f(x)+\frac{1}{2}\E\int_0^t|[\sigma_s(\sigma_s\sigma_s^\ast)^{-1}](X_s^x)(b^{(0)}_s(X_s^x)-b^{(0),n}_s(X_s^x))|^2\d s.
\end{align*}
So, Pinsker's inequality implies
\begin{align}\label{Pin}
\nonumber |P_t^{\<n\>} f(x)-P_tf(x)|^2&\leq \E\int_0^t|[\sigma(\sigma\sigma^\ast)^{-1}](X_s)(b^{(0)}_s(X_s)-b^{(0),n}_s(X_s))|^2\d s\\
& \leq \delta_2^{-1}\int_0^t\E|(b^{(0)}_s(X_s^x)-b^{(0),n}_s(X_s^x))|^2\d s
\end{align}
Applying \cite[Corollary 2.2.2]{BKR}, we conclude that there exists a density function $p_x(t,y)$ such that $\L_{X_t^x}(\d y)\d t=p_x(t,y)\d y\d t$, $t\in(0,T)$. By \eqref{Pin} and the dominated convergence theorem, we derive
\begin{align*}
\lim_{n\to\infty}|P_t^{\<n\>} f(x)-P_tf(x)|^2&\leq  \lim_{n\to\infty}\delta_2^{-1}\int_0^t\int_{\R^d}|(b^{(0)}_s(y)-b^{(0),n}_s(y))|^2p_x(s,y)\d y\d s=0.
\end{align*}
So, \eqref{kct} holds and the proof is completed.
\end{proof}
    Repeating the proof of \cite[Lemma S3]{HYY}, we have
\begin{lem}\label{yap} Assume {\bf(D)} and $b_t^{(0)}$ and $\tilde{b}_t^{(0)}$ are Lipschitz continuous with Lipschitz constant independent of $t$. Let $X_{s,t}^{m,x}$ and $\tilde{X}_{s,t}^{m,x}$ be defined in {\bf Step (ii)} of the proof of Lemma \ref{gracd1}. Then for any $t\geq s, x\in\R^d$, there exists a subsequence $\{m_n\}_{n\geq 1}$ such that $\P$-a.s. $\lim_{n\to\infty} |X_{s,t}^{m_n,x}-X_{s,t}^x|=0$ and $\lim_{n\to\infty} |\tilde{X}_{s,t}^{m_n,x}-\tilde{X}_{s,t}^x|=0.$
\end{lem}
\begin{proof}

By \eqref{abc}-\eqref{jg4dv7} and noting that $b_t^{(0)}$ and $\tilde{b}_t^{(0)}$ are Lipschitz continuous with Lipschitz constant independent of $t$, the proof is completed by repeating the proof of \cite[Lemma S3]{HYY}.
\end{proof}
\section*{}
{\bf Data Availability Statement} Data sharing not applicable to this article as no datasets were generated or
analysed during the current study.
\section*{Declarations}
{\bf Conflict of Interests} The authors declare that they have no conflict of interest.

\end{document}